\documentclass[11pt]{amsart}

\usepackage{amsmath}
\usepackage{amsfonts}
\usepackage{amssymb} 
\usepackage{amsthm} 
\usepackage{mathptmx}
\usepackage{mathrsfs}
\usepackage{latexsym}
\usepackage{times}
\usepackage[isolatin]{inputenc}

\DeclareMathAlphabet{\mathpzc}{OT1}{pzc}{m}{it}

\newtheorem{thm}{Theorem}[section]
\newtheorem{lem}[thm]{Lemma}
\newtheorem{prop}[thm]{Proposition}
\newtheorem{cor}[thm]{Corollary}

\theoremstyle{definition}
\newtheorem{defn}[thm]{Definition}

\theoremstyle{remark}
\newtheorem{rem}[thm]{Remark}

\newcommand{\osp}{\mathfrak{osp}}
\newcommand{\spk}{\mathfrak{sp}}

\newcommand{\slk}{\mathfrak{sl}}

\newcommand{\Sk}{\mathfrak{S}}
\newcommand{\g}{\mathfrak{g}}

\newcommand{\hk}{\mathfrak{h}}

\newcommand{\nk}{\mathfrak{n}}
\newcommand{\sk}{\mathfrak{s}}

\newcommand{\Xc}{\mathcal{X}}
\newcommand\Fc{\mathcal F}
\newcommand\Lc{\mathcal L}
\newcommand\Dc{\mathcal D}
\newcommand\Hc{\mathcal H}

\newcommand\Pc{\mathcal P}
\newcommand\Cc{\mathcal C}

\newcommand\Sc{\mathcal S}
\newcommand\Wc{\mathcal W}

\newcommand\Uc{\mathcal U}

\newcommand{\Wb}{\mathsf{W}}
\newcommand{\Sb}{\mathsf{S}}

\newcommand\RR{\mathbb R}
\newcommand\CC{\mathbb C}
\newcommand\NN{\mathbb N}
\newcommand\ZZ{\mathbb Z}

\newcommand\KK{\mathbb K}

\newcommand{\Tr}{\operatorname{Tr}}
\newcommand{\ze}{{\scriptscriptstyle{\overline{0}}}}
\newcommand{\un}{{\scriptscriptstyle{\overline{1}}}}
\newcommand\gO{\g_{\ze}}
\newcommand\gI{\g_{\un}}

\newcommand{\Spo}{\operatorname{Sp}}

\newcommand{\ad}{\operatorname{ad}}
\newcommand{\exx}{\operatorname{e}}

\renewcommand\dfrac{\displaystyle \frac}

\newcommand{\Dom}{\operatorname{Dom}}

\newcommand{\sym}{\operatorname{Sym}}
\newcommand{\Id}{\operatorname{Id}} \renewcommand\hat\widehat
\renewcommand\tilde\widetilde 
\newcommand{\spa}{\operatorname{span}}

\newcommand{\IW}{\operatorname{IW}}
\newcommand{\Str}{\operatorname{Str}}
\newcommand{\RStr}{\operatorname{RStr}}
\newcommand{\Strw}{\operatorname{Str}_\Wb}
\newcommand{\Strwb}{\operatorname{Str}_{\overline{\Wb}}} 

\newcommand{\codim}{\operatorname{codim}}

\newcommand\sta{{ \ \star \ }}
\newcommand\tst{{ \ \mathop{\star}\limits_t \ }}
\newcommand\Lcst{{\mathop{\scriptscriptstyle{\Lc}}\limits_t}}
\newcommand\tstu{{ \ \mathop{\star}\limits_1 \ }}
\newcommand\tsti{{\mathop{\star}\limits_t}}
\newcommand\rp{{ \ \mathop{\times}\limits_\rho \ }}

\newcommand\adpo{\ad_{\tt P}}

\newcommand{\dss}{\displaystyle}

\begin{document}

\title[Supertrace and superquadratic Lie structure on the Weyl
  algebra]{Supertrace and superquadratic Lie structure on the Weyl
  algebra, with applications to formal inverse Weyl transform}

\author{Georges Pinczon, Rosane Ushirobira}

\address{Institut de Mathématiques de Bourgogne, Université de
  Bourgogne, B.P. 47870, F-21078 Dijon Cedex, France}

\email{gpinczon, rosane@u-bourgogne.fr}



\date{\today}

\begin{abstract} 
Using the Moyal $\star$-product and orthosymplectic supersymmetry, we
construct a natural non trivial supertrace and an associated non
degenerate invariant supersymmetric bilinear form for the Lie
superalgebra structure of the Weyl algebra $\Wb$. We decompose adjoint
and twisted adjoint actions. We define a renormalized supertrace and a
formal inverse Weyl transform in a deformation quantization framework
and develop some examples.
\end{abstract}

\maketitle

\section{Introduction}

Since Killing, it is known that the existence of a non degenerate
invariant symmetric bilinear form is a crucial property for a Lie
algebra. Let us call such a Lie algebra a {\em quadratic} Lie
algebra. These Lie algebras are of interest, not only when they are
finite dimensional, but infinite dimensional as well \cite{Kac}. A
corresponding notion of a {\em superquadratic} Lie superalgebra
exists: namely, it is a Lie superalgebra endowed with a non degenerate
invariant supersymmetric bilinear form. Besides, invariant bilinear
forms are often constructed from Traces or Supertraces (though it is
not always the case), even in the infinite dimensional case.

The main result of the present paper is:

\medskip

{\bf{\em Theorem} 1 :} The Weyl algebra is a superquadratic Lie
superalgebra with (essentially unique) bilinear form derived from a
supertrace.

\medskip

To explain the origin of this statement, let us consider the Weyl
algebra $\Wb$ in $2n$ generators realized as the polynomial algebra in
$2n$ indeterminates endowed with the Moyal $\star$-product denoted by
$\star$. Briefly, in a deformation quantization framework (see
\cite{BFFLS,Stern}), we quantize the natural Poisson bracket in an
explicit way (Moyal $\star$-product). The main advantage of this
quantization is the fact of being invariant with respect to the
natural embedding of the Lie superalgebra $\osp(1,2n)$ in $\Wb$
(\cite{Fronsdal, Musson,GPU}), i.e. with respect to the
orthosymplectic supersymmetry. The conjunction of both arguments,
explicit Moyal product and supersymmetry, is indeed a powerful
machinery that allows to deduce algebraic properties of the Weyl
algebra, as we shall show in this paper.

The Weyl algebra $\Wb$ is a Lie algebra with bracket denoted by $[ . ,
  . ]_{\scriptscriptstyle{\Lc}}$, and since it is naturally
$\ZZ_2$-graded, $\Wb$ is a Lie superalgebra with bracket denoted by
$[.,.]$.

We define the supertrace on the Weyl algebra $\Wb$ as the usual
evaluation at 0: 
\[ \Str(F):=F(0), \ \forall \ F \in \Wb.\]

and a bilinear form $\kappa$ on $\Wb$ as:
\[ \kappa(F,G) := \Str(F \sta G), \ \forall \ F, G \in \Wb.\]

Then we prove:

\medskip

{\bf{\em Theorem} 2:} 

\begin{enumerate}

\item $\Str$ is a supertrace on $\Wb$, that is $\Str([F,G])= 0$, for
  all $F$, $G \in \Wb$ and one has $\ker(\Str) = [\Wb, \Wb]$.

\item The bilinear form $\kappa$ is supersymmetric, non degenerate and
  invariant under the adjoint representation of the Lie superalgebra
  $\Wb$.

\end{enumerate}

\medskip

Theorem 1 is a consequence of Theorem 2. Let us mention that obviously
there is no non trivial trace on the Lie algebra $\Wb$.

As another consequence of Theorem 2, we deduce a decomposition of
$\Wb$:
\[ \Wb = \KK \oplus [ \Wb, \Wb],\]

recovering in a natural way a nice result of I. M. Musson \cite
{Musson}.

Concerning the bilinear form $\kappa$, given the decomposition of
$\Wb$ into homogeneous factors (for the commutative product), $\Wb =
\bigoplus_{k \geq 0} \Sb^k$, we show that
\[\kappa(\Sb^k, \Sb^\ell) =\{0\}, \ k \neq \ell,\]

so that the restriction of $\kappa$ to each $\Sb^k$ is non degenerate,
and provides an explicit (orthogonal or symplectic, according to the
parity of $k$) invariant bilinear form for the standard simple action
of $\spk(2n)$ on $\Sb^k$.

We study in addition, the decomposition of the adjoint $\ad$ and
twisted adjoint $\ad'$ actions of the Lie superalgebra $\Wb$ on itself:

\medskip

{\bf{\em Theorem} 3:} 

\begin{enumerate}

\item Under the adjoint action, $\Wb$ decomposes as $\Wb = \KK \oplus
  [ \Wb, \Wb]$. Besides, $[\Wb, \Wb]$ is a simple $\ad(\Wb)$-module.

\item $\Wb$ is a simple $\ad'(\Wb)$-module and has a non
  degenerate invariant supersymmetric bilinear form.

\end{enumerate}

\medskip

The bilinear form associated to the $\ad'$-action is deduced from
$\kappa$, and extends the bilinear form used to construct the natural
embedding of $\osp(1,2n)$ in $\Wb$. As a corollary:

\medskip

{\bf{\em{Corollary:}}} $[\Wb, \Wb]$ is a simple superquadratic Lie
superalgebra.

\medskip

Let us quote that the simplicity of $[\Wb, \Wb]$ was known some time
ago, as a combination of a result by S. Montgomery \cite {Mont}
proving that $[\Wb, \Wb]/ \KK$ is simple and a result by I. M. Musson
\cite{Musson} proving that $\Wb = \KK \oplus [\Wb, \Wb]$. Our proof
using a supertrace and the Moyal product is direct and completely
different of the initial proofs in \cite{Mont} and \cite{Musson}.

Finally, we reinterpret the supertrace in a Renormalization Theory
context. In a few words, let $\Pc$ be the polynomial algebra in $n$
indeterminates: the Weyl algebra $\Wb$ acts on $\Pc$ as differential
operators with polynomial coefficients. It has been shown in
\cite{Pinczon1} that any linear operator on $\Pc$ is in fact a
differential operator, eventually of infinite order. A slight
improvement provides a very explicit remarkable formula:

\medskip

{\bf{\em Theorem} 4:} If $T \in \Lc(\Pc)$, then 
\[T = \sum_N \dfrac{1}{N!} \left(  m \circ (T \otimes \Sc) \circ \Delta
(x^N) \right) \ \frac{\partial^N}{\partial x^N}\]

where $m$ is the product, $\Delta$ the coproduct and $\Sc$ the
antipode of $\Pc$. 

\medskip

Let $\Str$ be the ordinary supertrace defined on the ideal of finite
rank operators $\Lc_f(\Pc)$ and $\Str_\Wb$ be the supertrace
previously defined on $\Wb$. We show the following theorem:

\medskip

{\bf{\em Theorem} 5:} Let $T = \sum_I \alpha_I(Q) \sta P^I$. If $T \in
\Lc_f(\Pc)$, then:
\[ \Str(T) = \frac{1}{2^n} \sum_I \Strw (\alpha_I(Q) \sta P^I).\]

\medskip

Therefore we can construct a natural extension of $\Str$ to
$\Lc_f(\Pc) \oplus \Wb$ (in fact, to a bigger subspace that will not
be discussed here). For this renormalized extension, that we note
$\RStr$, one has:
\[ \RStr(\Id) = \frac{1}{2^n},\]
a formula that renormalizes $\Str(\Id) = \infty - \infty$ obtained by
wildly applying the definition of $\Str$. Moreover, the power $n$ in
this formula is precisely the dimension of the underlying variety so
our renormalization of supertrace has a geometric flavor. At last, we
give some ideas of what a formal inverse Weyl transform could be, and compute
examples.

We have tried to make this paper as self-contained as possible. For
instance, we give a short introduction to the Moyal $\star$-product in
Section \ref{Section01}, and in Section \ref{Section02}, an
introduction to the embedding of $\osp(1,2n)$ in $\Wb$, together with
the decompositions of the corresponding adjoint and twisted adjoint
actions. Not only because all this material provides tools used all
along this paper, but also for the convenience of the reader, who is
not forced to be an expert in deformation quantization theory,
orthosymplectic supersymmetry, etc.

\subsection*{Remarks}

\begin{enumerate}

\item Theorem 2 results from Proposition \ref{1.24} and Theorem
  \ref{4.4}. Theorem 3 corresponds to Theorem \ref{3.6}. Theorem
  \ref{5.7.2} gives Theorem 4 and the proof of Theorem 5 is Theorem
  \ref{5.9}.
  
\item In Sections 1 to 4 of this paper, $\KK$ denotes a field of
  characteristic zero (not necessarily algebraically closed). In
  Section 5, $\KK = \RR$ or $\CC$.

\item There are so many references on deformation quantization that we
  cannot quote all of them. The reader should refer to the beautiful
  paper \cite{BFFLS} that is the beginning (and much more) of this
  theory and to \cite{Stern} for a complete history with 164
  references.

\end{enumerate}

\section{Moyal $\star$-products} \label{Section01}

Moyal $\star$-products are the first examples of Deformation
Quantization of Poisson Brackets (see e.g. \cite{BFFLS,Stern}). In this
Section, we recall their well-known properties, giving proofs for the
convenience of the reader.

Let $V$ be a vector space with a basis $\{ P_1, Q_1, \dots, P_n,
Q_n\}$ and $\{ p_1, q_1, \dots, p_n, q_n\}$ its dual basis. 

We denote by $\Sb$ the commutative algebra $\Sb:= \sym(V^*)$ with the
usual grading $\Sb = \oplus_{k \geq 0} \Sb^k$ and the Poisson bracket:
\[ \{F,G\} = \sum_{i=1}^n \left( \dfrac{\partial F}{\partial p_i}
\dfrac{\partial G}{\partial q_i} - \dfrac{\partial F}{\partial q_i}
\dfrac{\partial G}{\partial p_i} \right), \ \forall \ F, G \in \Sb.\]

Denoting by $\{ X_1, \dots, X_{2n} \}$ the given basis of $V$ and by
$\{ x_1, \dots, x_{2n} \}$ its dual basis, we introduce a duality that
identifies $\Sb^*$ and the commutative algebra of formal power series
$\Fc:= \KK[[X_1, \dots, X_{2n}]]$ by:
\[\langle x^I \mid X^J \rangle := \delta_{I,J} \ I! \]

where $I$, $J \in \NN^{2n}$, $x^I = x_1^{i_1} \dots x_{2n}^{i_{2n}}$,
$X^J = X_1^{j_1} \dots X_{2n}^{j_{2n}}$ and $I! = i_1! \dots i_{2n}
!$.

The following properties result from a straightforward verification:

\begin{prop} \label{1.3}
\hfill

\begin{enumerate}

\item For all $v \in V$ and $F \in \Sb$, $\langle F \mid \exx^v
  \rangle= F(v)$.

\item For all $v \in V$ and $F \in \Sb$, 
\[ \langle F \mid X^I \exx^v \rangle = \left\langle
\left. \dfrac{\partial^I F}{\partial x^I} \ \right| \ \exx^v
\right\rangle = \dfrac{\partial^I F}{\partial x^I} (v),\]

where $\dfrac{\partial^I F}{\partial x^I} := \dfrac{\partial^{i_1 +
    \dots i_{2n}} F}{\partial x_1^{i_1} \dots \partial
  x_{2n}^{i_{2n}}}$.

\end{enumerate}

\end{prop}

Notice that (1) is Taylor's formula and (2) means that the transpose
of the operator $\dfrac{\partial^I}{\partial x^I}$ of $\Sb$ is the
operator of multiplication by $X^I$ of $\Fc$.

Define an operator $\wp \colon \Sb \otimes \Sb \to \Sb \otimes \Sb$
by:
\begin{equation}\label{1.4}
\wp( F \otimes G) := \sum_{i=1}^n \left( \dfrac{\partial F}{\partial
  p_i} \otimes \dfrac{\partial G}{\partial q_i} - \dfrac{\partial
  F}{\partial q_i} \otimes \dfrac{\partial G}{\partial p_i} \right),
\ \forall \ F, G \in \Sb.
\end{equation}

Since $\Sb \otimes \Sb = \KK[p_1, q_1, \dots, p_n, q_n, p_1', q_1',
  \dots, p_n', q_n']$, we have:
\[ (\Sb \otimes \Sb)^* = \KK[[ P_1, Q_1, \dots, P_n, Q_n, P_1', Q_1',
  \dots, P_n', Q_n']]. \]
 
If $D = {}^t \wp$, from Proposition \ref{1.3} it follows:
\[D(A \otimes B) = \left( \sum_{i=1}^n (P_i \otimes Q_i - Q_i \otimes
P_i) \right) \cdot \ A \otimes B, \ \forall \ A, B \in \Fc, \]

so that $D$ is the operator of multiplication by $d(\Xc, \Xc') =
\sum_{i=1}^n (P_i Q_i' - Q_i P_i')$ of $(\Sb \otimes \Sb)^*$ where
$\Xc=(P_1, Q_1, \dots, P_n, Q_n)$ and $\Xc'=(P_1', Q_1', \dots, P_n',
Q_n')$.

Denote by $m$ the product of $\Sb$. We can now define:

\begin{defn} \label{1.6}
For $F$, $G \in \Sb$, the {\em Moyal $\star$-product} is defined as:
\[F \tst G := \left( m \circ \sum_{k \geq 0} \dfrac{t^k}{2^k \ k!}
\ \wp^k \right) (F \otimes G).\]
\end{defn}

Let $\Delta$ be the coproduct of $\Fc$ associated to $m$. With the
usual obvious abuse of notation, one can write:
\[\Delta(A) = A(\Xc + \Xc'), \forall \ A \in \Fc.\]

So if we denote by $m_\tsti$ the Moyal $\star$-product and by
$\Delta_\tsti$ the associated coproduct $\Delta_\tsti = {}^t m_\tsti$,
we get:
\[ \Delta_\tsti (A) = \sum_{k \geq 0} \dfrac{t^k}{2^k \ k!}
\ D^k (A(\Xc + \Xc')) = \exx^{\frac{t}{2} d(\Xc,\Xc')} \ A(\Xc +
\Xc'), \forall \ A \in \Fc.\]

\begin{prop} \hfill 

\begin{enumerate}
\item $\Delta_\tsti$ is coassociative.

\item  $m_\tsti$ is associative.

\end{enumerate}

\end{prop}

\begin{proof}
Let $A \in \Fc$. One has:
\begin{eqnarray*}
(\Id \circ \Delta_\tsti)(\Delta_\tsti(A)) &=& (\Id \circ
    \Delta_\tsti)\left( \exx^{\frac{t}{2} d(\Xc, \Xc')} \ A(\Xc +
    \Xc') \right) \\ &=& \exx^{\frac{t}{2} d(\Xc', \Xc'')}
    \exx^{\frac{t}{2} d(\Xc, \Xc'+ \Xc'')} \ A(\Xc + \Xc'+ \Xc'')
\end{eqnarray*}

\begin{eqnarray*}
(\Delta_\tsti \circ \Id )(\Delta_\tsti(A)) &=& \exx^{\frac{t}{2}
    d(\Xc, \Xc')} \Delta_\tsti(A) (\Xc + \Xc', \Xc'') \\ &=&
  \exx^{\frac{t}{2} d(\Xc, \Xc')} \exx^{\frac{t}{2} d(\Xc +\Xc',
    \Xc'')} \ A(\Xc + \Xc'+ \Xc'')
\end{eqnarray*}

and since $d$ is bilinear, one obtains $(\Id \circ \Delta_\tsti) \circ
\Delta_\tsti = (\Delta_\tsti\circ \Id ) \circ \Delta_\tsti$, so
$\Delta_\tsti$ is coassociative. From $\Delta_\tsti = {}^t m_\tsti$,
it results $m_\tsti$ associative.

\end{proof}

We introduce the notation:
\begin{equation}\label{1.10}
C_k(F,G) = \frac{1}{2^k k!} \left( m \circ \wp^k \right) (F \otimes
G), \ \forall \ F, G \in \Sb.
\end{equation}

Then we can write:
\[ F \tst G = FG + \dfrac{t}{2}  \{ F,G \} + t^2 C_2(F,G)+ \dots,
\ \forall \ F, G \in \Sb.\] 

Thus the Moyal $\star$-product is a deformation of the commutative
product of $\Sb$ and also a deformation quantization of the Poisson
bracket.

We remark that $C_k$ is a bidifferential operator of order
$(k,k)$, so:

\begin{prop} \label{1.12}
For all $F \in \Sb^f$ and $G \in \Sb^g$, $C_k(F,G) \in \Sb^{f+g-2k}$
(with $\Sb^\ell = \{0 \}$ when $\ell < 0$) and
\[ F \tst G = \sum_{k=0}^{\min(f,g)} \ t^k \ C_k(F,G).\]
\end{prop}

Consider $\tau$ the twist $\tau(F \otimes G) := G \otimes F$, $\forall
\ F,G \in \Sb$. Then $\wp \circ \tau = - \tau \circ \wp$, so $\wp^k
\circ \tau = (-1)^k \tau \circ \wp^k$ and one deduces a useful parity
property of the coefficients $C_k$:
\begin{equation}\label{1.13}
C_k(G,F) = (-1)^k C_k(F,G), \ \forall \ F,G \in \Sb.
\end{equation}

As a consequence, the Lie bracket associated to the Moyal
$\star$-product contains only odd terms:
\begin{equation}\label{1.14}
[F,G]_\Lcst = 2 \ \sum_{p \geq 0}  \ t^{2p+1} \ C_{2p+1} (F,G),
\ \forall \ F,G \in \Sb.
\end{equation}

Now, from Proposition \ref{1.12}, we obtain the following useful
property:

\begin{prop} \label{1.15}
Let $F \in \Sb$ of degree $\leq 2$. Then $[F,G]_\Lcst = t \{ F,G \}$,
for all $G \in \Sb$.
\end{prop}

Therefore $[p_i, q_i]_\Lcst = t$, for $i = 1, \dots, n$ and all other
brackets between $p_i$'s and $q_j$'s vanish.

Next we shall relate the Moyal $\star$-product to the Weyl algebra. For
this purpose, we need a Lemma:

\begin{lem}\label{1.16}
One has

\begin{enumerate}
\item For all $\varphi \in V^*$, $\varphi^{\tst k} =
  \underbrace{\varphi \tst \dots \tst \varphi}_{k \textrm{ times}} =
  \varphi^k$.

\item For all $\varphi_1, \dots, \varphi_k \in V^*$, $\varphi_1 \dots
  \varphi_k = \dfrac{1}{k!} \sum_{\sigma \in \Sk_k}
  \varphi_{\sigma(1)} \tst \dots \tst \varphi_{\sigma(k)}$.
\end{enumerate}

\end{lem}

\begin{proof}
To prove (1): by Proposition \ref{1.12}, one has $\varphi \tst \varphi
= \varphi^2+ \frac{t}{2} \{ \varphi, \varphi \}= \varphi^2$. Also
$\varphi \tst \varphi \tst \varphi = \varphi^3+ \frac{t}{2} \{
\varphi, \varphi^2 \} = \varphi^3$ and so on.

Then $(\lambda_1 \varphi + \dots + \lambda_k \varphi_k)^{\tst k} =
(\lambda_1 \varphi + \dots + \lambda_k \varphi_k)^k$, for all
$\lambda_1, \dots, \lambda_k \in \KK$. Identifying the coefficient of
the term $\lambda_1 \dots \lambda_k$ on each side, one obtains:
\[\sum_{\sigma \in \Sk_k} \varphi_{\sigma(1)} \tst \dots \tst
\varphi_{\sigma(k)} = \sum_{\sigma \in \Sk_k} \varphi_{\sigma(1)}
\dots \varphi_{\sigma(k)}\]

and (2) follows.
\end{proof}

\begin{rem} \label{1.16.1}
Consequently, since $\Sb^k$ is linearly generated by $\{ \varphi^k
\mid \varphi \in V^* \}$, we can conclude that $(\Sb, m_\tsti)$ is
generated by $V^*$ as an algebra.
\end{rem}

Let us denote by $\Wb(n)$, or $\Wb$ when there is no ambiguity, the
Weyl algebra generated by $\{p_1, q_1, \dots, p_n, q_n \}$ with
relations $p_i q_j-q_j p_i = \delta_{ij}$, $p_i p_j-p_j p_i = q_i
q_j-q_j q_i=0$, for $i, j =1, \dots, n$ (see \cite{Dixmier}). A
well-known fact is that $\Wb$ is the algebra of polynomial
differential operators on $\KK [y_1, \dots, y_n]$ with $p_i =
\dfrac{\partial}{\partial y_i}$ and $q_i = y_i
\scriptscriptstyle{\times} \cdot \ $. We can also realize $\Wb$ as the
quotient algebra $T(V^*)/ I$ where $I$ is the ideal in the tensor
algebra $T(V^*)$ generated by the above relations. Therefore we can
consider that $V^* \subset \Wb$ as will be done in the sequel.

Denote by $\rho \colon \Sb \to \Wb $ the {\em symmetrization map}
defined as:
\[ \rho(\varphi_1 \dots \varphi_k) = \dfrac{1}{k!} \sum_{\sigma \in
  \Sk_k} \varphi_{\sigma(1)} \varphi_{\sigma(2)} \dots
\varphi_{\sigma(k)}\]

where the right hand side products are computed in $\Wb$. Notice that
$\rho$ is well defined since it is the composition map of the usual
symmetrization map from $\Sb = \sym(V^*)$ into $T(V^*)$ with the
canonical map from $T(V^*)$ onto $\Wb = T(V^*)/I$.

Using the symmetrization $\rho$, we pull back the product of $\Wb$ on
$\Sb$ by:
\[F \rp G := \rho^{-1} (\rho(F) \ \rho(G)), \ \forall \ F, G \in
\Sb.\]

We fix the value $t=1$ of the parameter, and denote by $\star$ the
$\tstu$-Moyal product.

\begin{prop}\label{1.19}
One has $\rp = \star$. Therefore $F \sta G = \rho^{-1} (\rho(F)
\ \rho(G))$, for all $F, G \in \Sb$ and the algebra $(\Sb, \star)$ is
isomorphic to the Weyl algebra $\Wb$.
\end{prop}

\begin{proof}
It is clear from the definition that $\varphi^{\rp k} = \varphi^k$ for
all $\varphi \in V^*$, $k \in \NN$. As in the proof of Lemma
\ref{1.16}, one can write $\varphi_1 \dots \varphi_k = \frac{1}{k!}
\sum_{\sigma \in \Sk_k} \varphi_{\sigma(1)} \rp \dots \rp
\varphi_{\sigma(k)}$ and using the same Lemma once more, one has
$\varphi^{\rp k} = \varphi^{\star k}$, $\forall \varphi \in V^*$, $k
\in \NN$, and $\varphi_1 \rp \dots \rp \varphi_k = \varphi_1 \sta
\dots \sta \varphi_k$, $\forall \ \varphi_1, \dots, \varphi_k \in
V^*$. It results that $(\Sb, \rp)$ is generated by $V^*$ as an algebra
and the same holds for $(\Sb, \star)$ thanks to Remark
\ref{1.16.1}. So in order to prove that $\rp = \star$, we need only to
prove that they do coincide on the linear generators of $\Sb$, i.e. on
the monomials $\varphi_1 \rp \dots \rp \varphi_k = \varphi_1 \sta
\dots \sta \varphi_k$, for all $\varphi_1, \dots, \varphi_k \in
V^*$. But this is trivial.
\end{proof}

As a consequence of Proposition \ref{1.19}, we can identify the Weyl
algebra $\Wb$ and the Moyal algebra $\Sb$ endowed with the
$\star$-product. Such an identification reveals to be useful, since
the Moyal $\star$-product provides an explicit formula for the product
of the Weyl algebra, i.e. for the product of polynomial differential
operators. To illustrate, when $n=1$ one has the following explicit
formula:
\[ F \tst G = \sum_{k \geq 0} \frac{t^k}{2^kk!} \left(
\sum_{r+s=k} (-1)^s \binom{k}{s} \frac{\partial^k F}{\partial p_1^r
  \partial q_1^s} \frac{\partial^k G}{\partial p_1^s \partial q_1^r}
\right).\]

Consider now the basis $\{ q_1^i \sta p_1^j \mid (i,j) \in \NN^2 \}$
of $\Wb(1)$. By an easy computation, one obtains the following
expression in terms of classical orthogonal polynomials:
\begin{equation} \label{1.9.2}
q_1^i \sta p_1^j = \begin{cases} (-1)^j \ \dfrac{j!}{2^j}
  \ L_j^{(i-j)} (2p_1q_1) \ q_1^{i-j}, \text{ if } i \geq j
  \\ \\ (-1)^i \ \dfrac{i!}{2^i} \ L_i^{(j-i)} (2p_1q_1) \ p_1^{j-i},
  \text{ if } i \leq j \end{cases}.
\end{equation}

where $L_\beta^{(\alpha)}$ is the Laguerre polynomial (see
e.g. \cite{Szego}). Notice that this gives also the expression of the
basis $\{ q_1^{i_1} \sta p_1^{i_1} \sta \dots \sta q_n^{i_n} \sta
p_n^{i_n} \mid (i_1, \dots, i_n) \in \NN^n \}$ of $\Wb(n)$ since
$q_1^{i_1} \sta p_1^{i_1} \sta \dots \sta q_n^{i_n} \sta p_n^{i_n} =
q_1^{i_1} \sta p_1^{i_1} \cdot \ldots \cdot q_n^{i_n} \sta 
p_n^{i_n}$.

Let us now define a filtration and also a $\ZZ_2$-grading of
$\Wb$. For the filtration, we keep the filtration of $\Sb$, that is,
we set:
\[\Wb_k := \dss \bigoplus_{r \leq k} \ \Sb^r.\]

Using (\ref{1.10}) and Proposition \ref{1.12}, it is clear that this
is indeed a filtration of $\Wb$ and that the associated graded algebra
is $\Sb$.

What about the $\ZZ_2$-grading of $\Wb$? It can be defined in two
ways: first one, consider $\Wb$ as the algebra of polynomial
differential operators on $\KK[y_1, \dots, y_n]$. Define an element of
$\Wb$ to be {\em even} if it maps even polynomials into even
polynomials and odd polynomials into odd ones. An {\em odd} element
takes even polynomials into odd polynomials and vice versa. Evidently,
this defines a $\ZZ_2$-grading on $\Wb$ and one has
$\deg_{\ZZ_2}(p_1^{r_1} \sta q_1^{s_1} \sta \dots p_n^{r_n} \sta
q_n^{s_n}) \equiv (\sum_{i=1}^n r_i + \sum_{i=1}^n s_i) ($ mod $2 )$.

On the other hand, thanks to Proposition \ref{1.12}, one can define a
$\ZZ_2$-grading on $\Wb$ by
\[\Wb_\ze:= \bigoplus_{k \geq 0} \Sb^{2k} \text{ and }\Wb_\un :=
\bigoplus_{k \geq 0} \Sb^{2k+1}.\] 

This $\ZZ_2$-grading is exactly the preceding one. Therefore $\Wb$ can
be endowed with two Lie structures:

\begin{itemize}
\item a Lie algebra structure given by 
\[[F,G]_{\scriptscriptstyle{\Lc}}:= F \sta G - G \sta F, \text{ for
  all } F, G \in \Wb.\]

From (\ref{1.14}), one has $[\Wb_k, \Wb_\ell
]_{\scriptscriptstyle{\Lc}} \subset \Wb_{k + \ell-2}$, so $\Wb$ is a
filtered Lie algebra for the shifted filtration $\Wb_k^\Lc :=
\Wb_{k+2}$ and its associated graded Lie algebra is $\Sb$ endowed with
the Poisson bracket (using (\ref{1.14}) once again). A widely known
fact is that $[\Wb, \Wb]_{\scriptscriptstyle{\Lc}} = \Wb$ (we shall
give a proof in Section \ref{Section03}), so there does not exist a
non zero trace map on $\Wb$ satisfying $\Tr(F \sta G) = \Tr(G \sta
F)$, for $F$, $G \in \Wb$, i.e. $H^1(\Wb) = \{0\}$ for the Lie algebra
cohomology.

\item a Lie superalgebra structure given by 
\[ [F,G]:= F \sta G - (-1)^{fg} G \sta F, \text{ for all } F, G \in
\Wb.\]

with $f = \deg_{\ZZ_2}(F)$ and $g = \deg_{\ZZ_2}(G)$.

\end{itemize}

A significant difference between the Lie algebra and the Lie
superalgebra structures is set in the following Proposition:

\begin{prop} \label{1.24}
Let $\Str$ be the evaluation at 0 of the commutative algebra $\Sb$,
\[ \Str(F) := F(0), \ \text{ for all } \ F \in \Sb.\]

Then $\Str$ is a supertrace on $\Wb$, that is, $\Str$ is homogeneous
and satisfies
\[ \Str([F,G]) =0, \ \text{ for all } \ F, G \in \Wb.\]
\end{prop}

\begin{proof}
From Proposition \ref{1.12},
\[ F \sta G = \sum_{k=0}^{\min(f,g)} \ C_k(F,G)\]

with $C_k(F,G) \in \Sb^{f+g-2k}$ for $F \in \Sb^f$ and $G \in
\Sb^g$. One has $C_k(F,G) \in \KK$ if and only if $k =
\frac{f+g}{2}$. Since $k$ runs from 0 to $\min(f,g)$, one can deduce:
either $f \neq g$ and then $\Str(F \sta G) = 0$, or $f =g$ and $\Str(F
\sta G)= C_f(F,G)$. In this later case, if $f$ is even, then by
(\ref{1.13}), $C_f(F,G) = C_f(G,F)$, so $\Str(F \sta G) = \Str(G \sta
F)$. If $f$ is odd, $C_f(F,G) = - C_f(G,F)$ (\ref{1.13}), so $\Str(F
\sta G) = - \Str(G \sta F)$. Either way, $\Str([F,G]) = 0$. Finally,
$\Str(F) = 0$ if $F$ is odd, so $\Str$ is homogeneous.
\end{proof}

\begin{rem}
Assume that $\KK = \RR$ or $\CC$. The Moyal $\star$-product is clearly
not restricted to live only on polynomial functions: obviously,
assuming that $t$ is a formal parameter, Definition \ref{1.6} defines
an associative deformation of $\Cc^\infty (V)$. All above properties
are true, in the formal sense, thanks to the density of polynomials
inside $\Cc^\infty(V)$ (with its usual Fréchet topology, see
e.g. \cite{Treves}).

On the other hand, if $F$, $G \in \Sc(V)$ (fast decreasing smooth
functions), one has $\int F \sta G = \int F \ G$, so one can define a
trace by $\Tr(F) = \int F$ and this has very important consequences
(see e.g. \cite{Connes, Stern}).
\end{rem}

\section{Embedding $\spk(2n)$ and $\osp(1,2n)$ into
  $\Wb$}\label{Section02}

In the sequel, we denote by $\g$ the Lie superalgebra $\osp(1,2n)$ and
by $\g_\ze$ its even part, i.e. $\g_\ze = \spk(2n)$. The Weyl algebra
$\Wb$ is endowed with the super bracket $[F,G] = F \sta G - (-1)^{fg}
G \sta F$, for $F \in \Sb^f$, $G \in \Sb^g$. Denote by $\ad$ the
corresponding adjoint representation.

There exists a well-known embedding of $\g_\ze$ into $\Wb$ given by:

\begin{prop} \label{2.1}
Given $X \in \Sb^2$, denote by $\adpo(X):= \{X, . \}$. Then 
\[ \adpo (X)|_{\Sb^1} = \ad (X)|_{\Sb^1} \in \spk(2n).\]

Moreover $[\Sb^2, \Sb^2] = \{ \Sb^2, \Sb^2 \} \subset \Sb^2$, so that
$\Sb^2$ is a Lie algebra for the Poisson bracket or the
$\star$-bracket and the map 
\[X \mapsto \adpo(X)|_{\Sb^1}\]

is a Lie algebra isomorphism from $\Sb^2$ onto $\spk(2n)$.
\end{prop}

\begin{proof}
Set a non degenerated skew symmetric bilinear form $\Phi$ on $\Sb^1$
by $\Phi(\varphi,\psi) = \frac12 \{ \varphi, \psi \} = \frac12 [
  \varphi, \psi ]_{\scriptscriptstyle{\Lc}}$, for $\varphi$, $\psi \in
\Sb^1 = V^*$. Notice that when $X \in \Sb^2$, one has $\adpo(X) =
\ad(X)$ by Proposition \ref{1.15}. Then $[\Sb^2, \Sb^2] = \{ \Sb^2,
\Sb^2\} \subset \Sb^2$ which implies that $\Sb^2$ is a Lie algebra for
the Poisson bracket ( = $\star$-bracket). Since $\adpo(X)$ is a
derivation for the Poisson bracket, the bilinear form $\Phi$ is
invariant and that means $\adpo(X)|_{\Sb^1} \in \spk(2n)$. Let
$\theta$ be the map from $\Sb^2$ to $\spk(2n)$ defined as $\theta(X) =
\adpo(X)|_{\Sb^1}$. So $\theta$ is clearly a Lie algebra homomorphism,
which is injective since $\adpo(X)|_{\Sb^1} =0$ implies $X \in \KK$,
thus $X =0$. Besides $\theta$ is onto since $\dim(\Sb^2) =
\dim(\spk(2n))$.
\end{proof}

In the remaining of the paper, we shall identify $\g_\ze= \spk(2n)$
and $\Sb^2$ via the isomorphism in Proposition \ref{2.1}.

Let us define an embedding of $\g$ into the Weyl algebra $\Wb$
extending the previous embedding of $\g_\ze$ into $\Wb$. Long ago,
this embedding was already well-known among mathematical physicists
(e.g. \cite{Fronsdal}) and it was used for instance to develop
singleton Anti-de-Sitter theories \cite{Fronsdal}. This embedding was also
described in \cite{Musson} and \cite{GPU}. We give the proof in
\cite{GPU}:

\begin{defn} \label{2.2.0}
The {\em twisted adjoint action} of the Lie superalgebra $\Wb$ on
itself is defined as:
\[ \ad'(F)(G) := F \sta G - (-1)^{f(g+1)} G \sta F, \ \forall \ F
\in \Sb^f, G \in \Sb^g.\]
\end{defn}

Such a twisted action is typical of supersymmetry (see
\cite{Arnaudon,Gorelik}). Remark that if $F \in \Wb_\ze$, then
$\ad'(F) = \ad(F)$.

\begin{prop} \label{2.2}
Let $\sk = \Sb^1 \oplus \Sb^2$ with $\ZZ_2$-grading induced by
$\Wb$. Then $\sk$ is a subalgebra of the Lie superalgebra
$\Wb$. Moreover $\sk \simeq \osp(1,2n)$.
\end{prop}

\begin{proof}
If $\varphi$, $\psi \in \Sb^1$, one has $\varphi \sta \psi = \varphi
\psi + \frac12 \{ \varphi, \psi \}$, so $[\varphi, \psi]=2 \varphi
\psi \in \Sb^2$. If $\varphi \in \Sb^1$ and $X \in \Sb^2$, one has
$[\varphi,X] = \{ \varphi, X\}$ (by Proposition \ref{1.15}), and $\{
\varphi, X\} \in \Sb^1$. So $\sk$ is a subalgebra of the Lie
superalgebra $\Wb$.

Set $V = V_\ze \oplus V_\un$ where $V_\ze =\KK$ and $V_\un =
\Sb^1$. It is easy to check that $\ad'(\sk)(V) \subset V$. For
example, if $\varphi$, $\psi \in \Sb^1$, one has $\ad'(\varphi)(\psi)
= [\varphi, \psi]_{\scriptscriptstyle{\Lc}} = \{\varphi, \psi \} \in \KK$ and
$\ad'(\varphi)(1) = 2 \varphi$.

Define a bilinear form $\Theta$ on $V$ by $\Theta(\varphi, \psi) =
\frac12 \{ \varphi, \psi \}$ and $\Theta(1,1) = -1$. It is clear that
$\Theta$ is a supersymmetric form, i.e. $\Theta(A,B)= (-1)^{ab}
\Theta(B,A)$ for all $A \in V_a$ and $B \in V_b$. Straightforward
computations show that $\Theta$ is $\ad'$-invariant. Therefore we
obtain a Lie superalgebra homomorphism $\phi$ that maps $X$ to
$\ad'(X)|_V$ from $\sk$ into $\osp(1,2n)$. Remark that $\ad'(X)|_V =0$
for $X \in \Sb^1$ implies that $\ad'(X) (1) = 2 X =0$, so $\phi$ is
injective and since $\dim(\sk) = \dim(\osp(1,2n))$, $\phi$ is an
isomorphism.
\end{proof} 

In the remaining of the paper, we shall identify 
\[\g = \osp(1,2n) \qquad \text{ and } \qquad \sk = \Sb^1 \oplus \Sb^2\] 

via the isomorphism given in Proposition \ref{2.2}. Then 
\[\g = \Sb^1 \oplus \Sb^2,  \qquad \g_\un = \Sb^1 \qquad \text{  and }
\qquad \g_\ze = \Sb^2 = [\g_\un, \g_\un] = \spk(2n).\]

The root system of $\g$ is easily deduced: consider 
\[ H_i= -\frac12 \ [p_i,q_i] \text{ for } 1 \leq i \leq n.\]

Set $\hk = \spa \{ H_1, \dots, H_n\}$. It turns out that $\hk$ is a
Cartan subalgebra of $\gO$ and $\g$. Let $\{ \omega_1, \dots, \omega_n
\}$ be the dual basis. The root vectors and corresponding roots are:

\begin{itemize}

\item in $\gI$: \hfill

$p_i$ with root $\omega_i$ (positive), $q_i$ with root $-\omega_i$.

\item in $\gO$: \hfill

$[p_i,q_j]$ $(i \neq j)$ with root $\omega_i - \omega_j$ (positive if
  $i <j$),

$[p_i,p_j]$ with root $\omega_i + \omega_j$ (positive),

$[q_i,q_j]$ with root $-(\omega_i + \omega_j)$.

\end{itemize}

The fundamental root system is $\{ \omega_i - \omega_{i+1} \ (i =1,
\dots, n-1), \omega_n \}$ with corresponding root vectors $[p_i,
  q_{i+1}]$ $(i =1, \dots, n-1)$ and $p_n$. These vectors generate the
subalgebra $\nk^+$ of $\g$ with basis the positive root vectors. Any
simple finite dimensional $\g$-module has a {\em highest weight
  vector} $v$ of {\em weight} $\lambda \in \hk^*$ satisfying if $Hv =
\lambda(H) v $ for all $H \in \hk$ and $\nk^+ v =0$.

\begin{prop}\label{2.3}
The $\g$-module $\adpo(\gO)|_{\Sb^k}$ is a simple module with highest
weight vector $p_1^k$ of weight $k \omega_1$.
\end{prop}

\begin{proof}
See \cite{PS} or \cite {Musson} for a proof using Weyl's formula. We
provide here a short direct proof: it is clear that $\adpo(\gO)
(\Sb^k) \subset \Sb^k$. Now given $\varphi \in \Sb^1$, $\varphi \neq
0$, one can construct a Darboux basis with respect to $\Phi$ such that
$\varphi$ is the first basis vector. It results that there exists $A
\in \Spo(2n)$ such that $\varphi = A (p_1)$. Denoting by $A$ as well
the corresponding isomorphism of the commutative algebra $\Sb$, one
has $\varphi^k = A ( p_1^k)$. Since $\Sb^k = \spa \{\varphi^k \mid
\varphi \in \Sb^1 \}$, one has $\spa(\Spo(2n))(p_1^k) = \Sb^k$. One
has a representation of $\Spo(2n)$ in $\Sb^k$, $A \mapsto A|_{\Sb^k}$
with differential $\adpo(\gO)|_{\Sb^k}$. As a consequence, the
$\gO$-submodule of $\Sb^k$ generated by $p_1^k$ is $\Sb^k$
itself. Since $p_1^k$ is a highest weight vector of weight $k
\omega_1$, $\Sb^k$ is a simple $\gO$-module with highest weight $k
\omega_1$.
\end{proof}

The next Proposition was proved in \cite{Pinczon2} for $\osp(1,2)$ and
in \cite{Musson} for $\osp(1,2n)$. We give a simplified proof that uses
the Moyal $\star$-product.

\begin{prop}\label{2.4}
Let $A_k = \Sb^{2k-1} \oplus \Sb^{2k}$, $k \geq 1$ and $A_0=\{ 0
\}$. Then $A_k$ is stable by $\ad(\g)$. Moreover, $A_k$ is a simple
$\g$-module with highest weight vector $p_1^{2k}$ of weight $2 k
\omega_1$.
\end{prop} 

\begin{proof}
Let $X \in \gI$ and $F \in \Sb^{2k}$. Then $\ad(X)(F)=\{X,F\}$ by
Proposition \ref{1.15}, hence $\ad(X)(F) \in \Sb^{2k-1}$. Now assume
that $F \in \Sb^{2k-1}$, then $\ad(X)(F) = X \sta F + F \sta X = 2 X F
\in \Sb^{2k}$. So $A_k$ is $\ad(\g)$-stable. By Proposition \ref{2.3},
$A_k = \Sb^{2k-1} \oplus \Sb^{2k}$ is its decomposition into isotypic
components under the action of $\ad(\gO)$. Any $\g$-submodule $U$ of
$A_k$ must be decomposed as $U = (U \cap \Sb^{2k}) \oplus (U \cap
\Sb^{2k-1})$ and since $\Sb^{2k}$ and $\Sb^{2k-1}$ are simple
$\gO$-modules, one of them is contained in $U$ if $U \neq \{0\}$ from
what one deduces that $U = A_k$ using the beginning of the
proof. Moreover, $p_1^{2k}$ is clearly a highest vector of weight $2 k
\omega_1$.
\end{proof}

Next we decompose the twisted adjoint representation of $\g$ in $\Wb$:

\begin{prop}\label{2.6}
Let $B_k = \Sb^{2k} \oplus \Sb^{2k+1}$, $k \geq 0$. Then $B_k$ is
stable by $\ad'(\g)$. Moreover, $B_k$ is a simple $\g$-module with
highest vector $p_1^{2k+1}$ of weight $(2k+1) \omega_1$.
\end{prop}

\begin{proof}
Let $X \in \gI$ and $F \in \Sb^{2k}$. Then $\ad'(X)(F) = 2 X F \in
\Sb^{2k}$ and if $F \in \Sb^{2k+1}$, then $\ad'(X)(F) = \{X,F\} \in
\Sb^{2k}$, so $B_k$ is $\ad'(\g)$-stable. By the same arguments in the
proof of the previous Proposition, one obtains easily that $B_k$ is a
simple $\g$-module and $p_1^{2k+1}$ is clearly a highest weight vector
of weight $(2k+1) \omega_1$.
\end{proof}

\begin{rem}
From Propositions \ref{2.4} and \ref{2.6}, it results that the only
homomorphism of $\g$-modules from $(\Wb, \ad(\g))$ into $(\Wb,
\ad'(\g))$ (or the other way around) is zero.
\end{rem}

\section{Decomposition of adjoint and twisted adjoint
  $\Wb$-modules} \label{Section03}

Let us recall our conventions: the Weyl algebra $\Wb$ is a Lie algebra
with bracket denoted by
\[ [A,B]_{\scriptscriptstyle{\Lc}} = A \sta B - B \sta A,  \forall
\, A, B \in \Wb.\]

Denote by $\ad_{\scriptscriptstyle{\Lc}}$ the corresponding adjoint
representation. Thanks to the $\ZZ_2$-gradation of the associative
algebra $\Wb$, there is a twisted adjoint representation of the Lie
algebra $\Wb$ defined as:
\[ \ad_{\scriptscriptstyle{\Lc}}'(A)(B)  = A \sta B - (-1)^a B \sta A,
\forall \, A \in \Wb_a, B \in \Wb.\]

Note that $\ad_{\scriptscriptstyle{\Lc}}'(\Wb_\ze)(\Wb_i) \subset
\Wb_i \ $ ($i=\overline{0}$, $\overline{1}$) and
$\ad_{\scriptscriptstyle{\Lc}}'(\Wb_\un)$ takes $\Wb_\ze$ into
$\Wb_\un$ and vice versa.

The Weyl algebra is also a Lie superalgebra with bracket denoted by
\[ [A,B] = A \sta B - (-1)^{ab} B \sta A, \forall \, A \in \Wb_a, B
  \in \Wb_b.\]

Denote by $\ad$ the corresponding adjoint representation. There is a
twisted adjoint representation of the Lie superalgebra $\Wb$ defined
as:
\[ \ad'(A)(B)  = A \sta B - (-1)^{a(b+1)} B \sta A, \forall \, A
\in \Wb_a, B \in \Wb_b.\]

This twisted adjoint representation was used in Proposition \ref{2.3}
to prove the embedding of $\g = \osp(1,2n)$ in $\Wb$.

The decomposition of $\ad_{\scriptscriptstyle{\Lc}}(\gO) =
\ad_{\scriptscriptstyle{\Lc}}'(\gO)$ is given in Proposition
\ref{2.3}, the decomposition of $\ad(\g)$ is given in Proposition
\ref{2.4} and in Proposition \ref{2.6} for $\ad'(\g)$.

We will now examine $\ad_{\scriptscriptstyle{\Lc}}$,
$\ad_{\scriptscriptstyle{\Lc}}'$ as Lie algebra representations of
$\Wb$, $\ad$ and $\ad'$ as Lie superalgebra representations of
$\Wb$. The main technical argument is given by the Theorem:

\begin{thm} \label{3.3}
Consider the $\ad= \ad_{\scriptscriptstyle{\Lc}}$-action of $\gO$ on
$\Wb$. One has:

\begin{enumerate}

\item The map $F \otimes G \mapsto F \sta G$ is a homomorphism of
  $\gO$-modules from $\Sb^\ell \otimes \Sb^m$ onto
  $\bigoplus_{k=0}^{\min(\ell, m)} \ \Sb^{\ell +m -2k}$, so one has
\[\Sb^\ell \sta \Sb^m =  \bigoplus_{k=0}^{\min(\ell, m)} \ \Sb^{\ell +m -2k}.\]

\item The map $F \otimes G \mapsto [F, G]_{\scriptscriptstyle{\Lc}}$
  is a homomorphism of $\gO$-modules from $\Sb^\ell \otimes \Sb^m$
  onto $\bigoplus_{2p+1 \leq \min(\ell, m)} \ \Sb^{\ell +m -2(2p+1)}$,
  so one has
\[ [\Sb^\ell, \Sb^m]_{\scriptscriptstyle{\Lc}} =\bigoplus_{2p+1 \leq
  \min(\ell, m)} \ \Sb^{\ell +m -2(2p+1)}.\]

\item The map $F \otimes G \mapsto [F, G]$ is a homomorphism of
  $\gO$-modules from $\Sb^\ell \otimes \Sb^m$ onto

\begin{itemize}

\item $\bigoplus_{2p+1 \leq \min(\ell, m)} \ \Sb^{\ell +m -2(2p+1)}$
  if $\overline{\ell} \ \overline{m} \equiv \overline{0}$, so in this
  case one has
\[ [\Sb^\ell, \Sb^m] =\bigoplus_{2p+1 \leq \min(\ell, m)} \ \Sb^{\ell
    +m -2(2p+1)}.\]

\item $\bigoplus_{2p+1 \leq \min(\ell, m)} \ \Sb^{\ell +m -4p}$ if
  $\overline{\ell} \ \overline{m} \equiv \overline{1}$, so in this
  case one has
\[ [\Sb^\ell, \Sb^m] =\bigoplus_{2p \leq \min(\ell, m)} \ \Sb^{\ell
    +m -4p}.\]

\end{itemize}

\item The map $C_k$ is a homomorphism of $\gO$-modules from $\Sb^\ell
  \otimes \Sb^m$ into $\Sb^{\ell +m -2k}$, so one has
\[ \begin{cases} C_k(\Sb^\ell, \Sb^m) = \Sb^{\ell +m -2k},
\text{ if } \ 0 \leq k \leq \min(\ell,m)\\ 0, \text{ otherwise. }
  \end{cases} \]

\end{enumerate}
\end{thm}

To prove the above Theorem, we need the following Lemma:

\begin{lem} \label{3.4}
Let $\hk$ be a Lie algebra, $\Uc =\Uc(\hk)$ its enveloping algebra and
$V$ a semisimple $\hk$-module which decomposes as $V = \oplus_{s=1}^r
V_s$ with $V_s$ simple and two by two non isomorphic. Let $v =
\sum_{s=1}^r v_s$ with $v_s \in V_s$ and $v_s \neq0$. Then $V =
\Uc(\hk) v$.
\end{lem}

\begin{proof}
The $\Uc$-module $V$ is semisimple and $V = \oplus_{s=1}^r V_s$ is its
decomposition into isotypic components. Let $\pi_s$ be the projection
of $V$ into $V_s$. Then $\pi_s$ is an element of the bicommutant of
the $\Uc$-module $V$ and by the Jacobson density theorem, there exists
$u \in \Uc$ such that $v_s = \pi_s(v)= u v$ and the result follows
from the simplicity of each $V_s$.
\end{proof}

\begin{proof} (Proof of Theorem \ref{3.3}) \hfill

For all $X \in \gO$, $\gO$ acts by $\ad(X)$, which is a derivation of
the $\star$-product. One has $F \sta G = \sum_{k=0}^{\min(\ell, m)}
C_k(F,G)$, for all $F \in \Sb^\ell$, $G \in \Sb^m$ with $C_k(F,G) \in
\Sb^{\ell + m -2k}$ which is irreducible under the
$\ad(\gO)$-action. It results that the map $F \otimes G \mapsto F \sta
G$ is a homomorphism of $\gO$-modules and so does any of the maps
$C_k$.

If $L = \bigoplus_{k=0}^{\min(\ell,m)} \Sb^{\ell+m-2k}$, denote by
$\rho$ the map $\rho \colon \Sb^\ell \otimes \Sb^m \to L$, $\rho(F
\otimes G) = F \sta G$. We compute $\rho(p_1^\ell \otimes q_1^m)$:
using (\ref{1.4}), one has $\wp(p_1^\ell \otimes q_1^m) = \ell m
p_1^{\ell-1} \otimes q_1^{m-1}$, so
\[\wp^k(p_1^\ell \otimes q_1^m) = \ell (\ell -1) \dots (\ell
-k+1)m(m-1)(m-k+1) p_1^{\ell-k} \otimes q_1^{m-k},\]

if $k \leq \min(\ell, m)$ and 0 otherwise. Then 
\[C_k(p_1^\ell,  q_1^m) = \frac{k!}{2^k} \binom{\ell}{k} \binom{m}{k}
p_1^{\ell-k} q_1^{m-k},\]

if $k \leq \min(\ell,m)$ and 0 otherwise. Thus
\[\rho(p_1^\ell \otimes q_1^m) = \sum_{k=0}^{\min(\ell,m)}
\frac{k!}{2^k} \binom{\ell}{k} \binom{m}{k} p_1^{\ell-k} q_1^{m-k}.\]

We can apply Lemma \ref{3.4} to the $\gO$-submodule $\rho(\Sb^\ell
\otimes \Sb^m)$ of $L$ with $v = \rho(p_1^\ell \otimes q_1^m)$ to
obtain that $\rho(\Sb^\ell \otimes \Sb^m) = L$.

But $C_k$ is a homomorphism of $\gO$-modules from $\Sb^\ell \otimes
\Sb^m$ into $\Sb^{\ell +m -2k}$ and one has $C_k(\Sb^\ell, \Sb^m) =
\{0\}$ if $k > \min(\ell, m)$. When $k \leq \min(\ell, m)$,
$C_k(\Sb^\ell, \Sb^m)$ is a non zero submodule of the simple module
$\Sb^{\ell +m-2k}$ since $C_k(p_1^\ell, q_1^m) \neq 0$. Therefore
$C_k(\Sb^\ell, \Sb^m) = \Sb^{\ell +m -2k}$.

The proof of (2) is completely similar, using (\ref{1.14}). The same
reasons can also be used to prove (3): the definition of $[.,.]$
implies that $[F,G] = 2 \sum_{p \geq 0} C_{2p+1} (F,G)$ if
$\overline{f} \overline{g} \equiv \overline{0}$ and $[F,G] = 2 \sum_{p
  \geq 0} C_{2p} (F,G)$ if $\overline{f} \overline{g} \equiv
\overline{1}$.
\end{proof}

\begin{rem} \hfill

\begin{enumerate}

\item Let us consider the case $n=1$, $\Wb = \Wb(1)$. Then $\gO =
  \slk(2)$ and $\Sb^\ell$ is the simple $(\ell +1)$-dimensional
  $\gO$-module, denoted by $D\left(\frac{\ell}{2} \right)$. Then
  $\dim(\Sb^\ell \otimes \Sb^m) = \dim \left(
  \bigoplus_{k=0}^{min(\ell,m)} \Sb^{\ell +m -2k} \right)$. So the map
  $F \otimes G \mapsto F \sta G$ of Proposition \ref{3.3}(1) is an
  isomorphism from $D \left( \frac{\ell}{2} \right) \otimes D \left(
  \frac{m}{2} \right)$ onto
\[ D \left( \frac{\ell}{2} +\frac{m}{2} \right) \oplus D \left(
\frac{\ell}{2} +\frac{m}{2} -1 \right) \oplus \dots \oplus D \left(
\left| \frac{\ell}{2} -\frac{m}{2} \right| \right),\]

providing an explicit and very handy formula for the computation of
Clebsch-Gordan coefficients. Such a formula was used for instance in
\cite{Arnal} to compute the commutation rules of high dimensional
Lie algebras when usual Clebsch-Gordan formulas were
hopeless. Unfortunately, it is easy to check that the $\star$-product
will provide only a partial decomposition of the tensor product
$\Sb^\ell \otimes \Sb^m$ when $n \geq 2$.

\item Using Proposition \ref{3.3}, one has the following identities:
\begin{eqnarray}\label{3.5.2}
& [\Sb^1, \Sb^k]_{\scriptscriptstyle{\Lc}} = \Sb^{k-1}, \quad [\Sb^2,
    \Sb^k]_{\scriptscriptstyle{\Lc}} = \Sb^{k} \ (k \geq 1), \\ \notag
  & [\Sb^3, \Sb^k]_{\scriptscriptstyle{\Lc}} = \Sb^{k-3} \oplus
  \Sb^{k+1} \ (k \geq 3), \\ \notag & [\Sb^3,
    \Sb^2]_{\scriptscriptstyle{\Lc}} = \Sb^3 \quad \text{ and } \quad
     [\Sb^3, \Sb^1]_{\scriptscriptstyle{\Lc}} = \Sb^2.
\end{eqnarray}
\begin{eqnarray}\label{3.5.3}
&[\Sb^1, \Sb^{2k}] = \Sb^{2k-1} \ (k \geq 1), \quad [\Sb^1,
    \Sb^{2k+1}] = \Sb^{2k+2} \ (k \geq 0), \\ \notag & [\Sb^3,
    \Sb^{2k-1}] = \Sb^{2k-2} \oplus \Sb^{2k+2} \ (k \geq 2), \quad
  [\Sb^3, \Sb^1] = \Sb^4,\\ \notag & [\Sb^3, \Sb^{2k}] = \Sb^{2k-3}
  \oplus \Sb^{2k+1} \ (k \geq 2) \quad \text{ and } \quad [\Sb^3,
    \Sb^2] = \Sb^3.
\end{eqnarray}

These identities turn out to be quite useful.
\end{enumerate}

\end{rem}

\begin{prop} \label{3.7}
One has 
\[\ker(\Str) = [ \gO, \Wb] = [\Wb, \Wb] = \bigoplus_{k \geq 1} \Sb^k.\]
\end{prop}

\begin{proof}
We need the following obvious result: 

{\em If $\hk$ is a Lie algebra and $U$ a non trivial simple
  $\hk$-module, then $\hk \ U = U$}.

Applying it to $\gO$ and the simple $\gO$-module $\Sb^k$ ($k \geq 1$)
we get that $[\gO, \Sb^k] = \Sb^k$ (alternatively, $[\Sb^2, \Sb^k] =
\Sb^k$ by Formula (\ref{3.5.2})). As a consequence, $[ \gO,
  \bigoplus_{k \geq 1} \Sb^k ] = \oplus_{k \geq 1} \Sb^k$. By
Proposition \ref{1.24}, one obtains $[\Wb, \Wb] \subset \ker(\Str)$,
so finally we conclude that $\bigoplus_{k \geq1} \Sb^k \subset [\Wb,
  \Wb] \subset \ker(\Str)$, but since $\codim \left( \bigoplus_{k \geq
  1} \Sb^k \right) = \codim( \ker(\Str)) = 1$, the result is proved.
\end{proof}

\begin{cor} \cite{Musson} 
One has $\Wb=\KK \oplus [\Wb,\Wb]$.
\end{cor}

\begin{rem}
The use of the supertrace $\Str$ is quite natural and provides a real
simplification of Musson's proof which does not use $\Str$. Moreover,
it provides an enlightening on the origin of the result.
\end{rem}

\begin{thm} \label{3.6} \hfill

\begin{enumerate}

\item The representation $(\ad_{\scriptscriptstyle{\Lc}}, \Wb)$ of the Lie algebra $\Wb$ is
  indecomposable with Jordan-Hölder series $\{0\} \subset \KK \subset
  \Wb$ and one has $\ad_{\scriptscriptstyle{\Lc}}(\Wb)(\Wb) = \Wb$.

\medskip 

\item The representation $(\ad_{\scriptscriptstyle{\Lc}}', \Wb)$ of the Lie algebra $\Wb$ is
  indecomposable, $[\Wb, \Wb]$ is a simple subrepresentation, $\Wb/
  [\Wb, \Wb]$ is the trivial representation and there exists a non
  trivial $[\Wb, \Wb]$-valued cocycle $\xi$ defined by $\xi(F) = 0$ if
  $F$ is even, and $\xi(F) = 2 F$ if $F$ is odd.

\medskip 

\item The representation $(\ad, [\Wb,\Wb])$ of the Lie superalgebra
  $\Wb$ is simple. Moreover, the $\ad(\Wb)$-module $\Wb$ decomposes as
  $\Wb = \KK \oplus [\Wb, \Wb]$.

\medskip 

\item The representation $(\ad', \Wb)$ of the Lie superalgebra $\Wb$
  is simple.

\end{enumerate}

\end{thm}

\begin{proof} \hfill

\begin{enumerate}

\item Let $M$ be a non trivial
  $\ad_{\scriptscriptstyle{\Lc}}(\Wb)$-module. Then $M$ can be
  decomposed into isotypic components under the
  $\ad_{\scriptscriptstyle{\Lc}}(\gO)$-action. By Proposition
  \ref{2.3}, there exists $k_0$ such that $\Sb^{k_0} \subset M$. If
  $k_0 > 0$, using (\ref{3.5.2}) one has $[\Sb^1,
    \Sb^{k_0}]_{\scriptscriptstyle{\Lc}} = \Sb^{k_0-1} \subset M$
  implying that $\bigoplus_{k \leq k_0} \Sb^k \subset M$. Using
  $\ad_{\scriptscriptstyle{\Lc}}(\Sb^3)$ and (\ref{3.5.2}), one
  deduces that $\Sb^{k_0+1} \subset M$ and repeating the argument,
  that $M = \Wb$.

If the only $k$ such that $\Sb^k \subset M$ is 0, then $M = \Sb^0 =
\KK$. Therefore, there is exactly one non trivial invariant subspace
namely $\KK$. It results that $\Wb / \KK$ is simple, that the
representation $(\ad_{\scriptscriptstyle{\Lc}}, \Wb)$ is
indecomposable, with Jordan-Hölder series $\{0\} \subset \KK \subset
\Wb$. To finish, notice that $\ad_{\scriptscriptstyle{\Lc}}(\Wb)(\Wb)$
is invariant, contains strictly $\KK$, so
$\ad_{\scriptscriptstyle{\Lc}}(\Wb)(\Wb) = \Wb$.

\medskip 

\item We start by proving that
  $\Str(\ad_{\scriptscriptstyle{\Lc}}'(F)(G)) = 0$, $\forall \ F$, $G
  \in \Wb$. If $F$ is even, or if $F$ and $G$ are odd, then
  $\ad'_{\scriptscriptstyle{\Lc}}(F)(G) = \ad(F)(G)$, so we apply
  Proposition \ref{1.24}. Now assume that $F \in \Sb^{k}$, $k$ odd and
  $G \in \Sb^\ell$, $\ell$ even. Then
  $\ad_{\scriptscriptstyle{\Lc}}'(F)(G) = F\sta G + G \sta F = 2
  \sum_{2s \leq \min(k, \ell)} C_{2s}(F,G)$ and $C_{2s}(F,G) \in
  \linebreak \Sb^{k+\ell -4s}$. But $k+\ell -4s > 0$, so
  $\ad'_{\scriptscriptstyle{\Lc}}(F)(G) \in \bigoplus_{r > 0} \Sb^r =
  \ker(\Str)$.

It results that $\ad_{\scriptscriptstyle{\Lc}}'(\Wb) (\Wb) \subset
\ker(\Str) = [\Wb, \Wb]$ and a fortiori, $[\Wb, \Wb]$ is stable. The
quotient $\Wb/ [\Wb, \Wb]$ is the trivial representation. If $F$ is
invariant under $\ad'_{\scriptscriptstyle{\Lc}}(\Wb)$, it is invariant
under $\ad_{\scriptscriptstyle{\Lc}}(\gO)$ and has to be a
constant. But $\ad'_{\scriptscriptstyle{\Lc}}(G)(1) = 2 G$, if $G$ is
odd, so $F =0$. Therefore the extension
\[ 0 \to [\Wb, \Wb] \to \Wb \to  \KK = \Wb / [\Wb, \Wb] \to 0\]

is non trivial with corresponding non trivial $[\Wb, \Wb]$-valued
cocycle $\xi$ defined by:
\[ \xi(F) = \ad'_{\scriptscriptstyle{\Lc}}(F)(1) = \begin{cases} 0, &
  \text{ if } F \text{ is even.} \\ 2F, & \text{ if } F \text{ is
    odd.}\end{cases}.\]

Assume now that $M$ is a non zero
$\ad'_{\scriptscriptstyle{\Lc}}(\Wb)$-invariant subspace. As in (1),
there exists there exists $k_0$ such that $\Sb^{k_0} \subset M$. Let
$X \in \Sb^1$ and $F \in \Sb^{k_0}$. One has
$\ad_{\scriptscriptstyle{\Lc}}'(X)(F) = 2 X F$, so
$\ad_{\scriptscriptstyle{\Lc}}'(\Sb^1)(\Sb^{k_0}) = \Sb^{k_0 +1}$ and
$\bigoplus_{k \geq k_0} \Sb^k \subset M$. Take $X \in \Sb^3$, then
$\ad_{\scriptscriptstyle{\Lc}}'(X)(F) = 2 X F + 2 C_2(X,F)$. As in the
proof of Theorem \ref{3.4}, it results that
$\ad_{\scriptscriptstyle{\Lc}}'(\Sb^3)(\Sb^{k_0}) = \Sb^{k_0 +3}
\oplus \Sb^{k_0 -1}$ if $k_0 > 1$, so $\Sb^{k_0 -1}$ is contained in
$M$ and then $\bigoplus_{k \geq 1} \Sb^k \subset M$. Thus, either $ M
= \bigoplus_{k \geq 1} \Sb^k = [\Wb, \Wb]$ or $M = \Wb$.

\medskip 

\item Thanks to Proposition \ref{3.7}, it remains to show that the
  representation \linebreak $(\ad, [\Wb,\Wb])$ is simple. This is an
  easy consequence of (\ref{3.5.3}): let $M$ be a non zero invariant
  subspace in $[\Wb, \Wb]$. Since $M$ is $\ad(\g)$-stable, it
  decomposes into the isotypic components of $[\Wb, \Wb] = \bigoplus_{k
    \geq 1} A_k$ with $A_k = \Sb^{2k-1} \oplus \Sb^{2k}$ (see
  Proposition \ref{2.4}). Hence $M = \underset{A_k \subset
    M}{\bigoplus_k} A_k$. Take $k_0$ to be the smallest $k$ in this
  decomposition. There are two cases:

\begin{itemize}

\item if $k_0 =1$, then $\g = A_1 \subset M$ and $[\gO, \bigoplus_{k
    \geq 1} \Sb^k]= \bigoplus_{k \geq 1} \Sb^k \subset M$, one has $M
  = [\Wb,\Wb]$ by Proposition \ref{3.7}.

\item if $k_0 > 1$, then $[\Sb^3, A_{k_0}] \subset M$. By
  (\ref{3.5.3}), we deduce that $\Sb^{2k_0 -3}$,  $\Sb^{2k_0 -2}$,
  $\Sb^{2k_0 +1}$ and  $\Sb^{2k_0 +2}$ are contained in $M$. But it
  results that $A_{k_0 -1} = \Sb^{2k_0 -3} \oplus \Sb^{2k_0 -2}$ is
  contained in $M$, a contradiction.

\end{itemize}

\medskip 

\item The proof of (4) is completely similar: one uses the
  decomposition of the representation $(\ad, \Wb)$ into isotypic
  components, $\Wb = \bigoplus_{k \geq0} B_k$ where $B_k = \Sb^{2k}
  \oplus \Sb^{2k+1}$ (see Proposition \ref{2.6}).

\end{enumerate}

\end{proof}

\section{Invariant bilinear forms for the adjoint and twisted adjoint
  actions of $\Wb$} \label{Section04} 

There is a bilinear form $\kappa$ canonically associated to the
supertrace on the Weyl algebra $\Wb$, namely:
\[ \kappa(F,G) := \Str(F \sta G), \ \forall \ F, G \in \Wb.\]

By Proposition \ref{1.24}, $\kappa$ is supersymmetric:
\[ \kappa(G,F)= (-1)^{fg} \kappa(F,G), \ \forall \ F \in \Sb^f, G \in
\Sb^g.\] 

and from its very definition, $\kappa$ is invariant under the adjoint
representation:
\[ \kappa(\ad(F)(G), H ) + (-1)^{fg} \kappa( G, \ad(F)(H)) = 0,  \ \forall
\ F \in \Sb^f, G \in \Sb^g.\]

Now, $\kappa$ will be really interesting if it is non degenerate, and
indeed this is the case:

\begin{thm}\label{4.4}
The bilinear form $\kappa$ on the Lie superalgebra $\Wb$ is non
degenerate supersymmetric and $\ad$-invariant. Moreover,
\[ \kappa( \Sb^\ell, \Sb^m) = \{0\}, \text{ if } \ell \neq m,\]
and the restriction of $\kappa$ to $\Sb^\ell$ is non degenerate. If
$\kappa'$ is another $\ad$-invariant bilinear form on $\Wb$, there exists
$\alpha \in \KK$ such that $\kappa' = \alpha \ \kappa$ on $[\Wb, \Wb]$.
\end{thm}

\begin{proof}
Given $F \in \Sb^\ell, G \in \Sb^m$ with $\ell < m$, then $\kappa
(F,G) = \Str \left( \sum_{k=0}^\ell C_k(F,G) \right)$ by Proposition
\ref{1.12}, with $C_k(F,G) \in \Sb^{\ell +m-2k}$ and since $\ell + m
-2k >0$, $\kappa(F,G) = 0$. Indeed, $\kappa(F,G) = 0$ if
$\overline{\ell} + \overline{m} = \overline{1}$.

To prove that $\kappa$ is non degenerate, we have to show it on
each component $\Sb^\ell$. Using Proposition \ref{1.12}, if $F$ $G \in
\Sb^\ell$, then $\kappa(F,G) = C_\ell (F,G)$. Define $\phi \colon
\Sb^\ell \to (\Sb^\ell)^*$ as $\phi(F)(G) = \kappa(F,G)$, for all $F$,
$G \in \Sb$. Since $\kappa$ is $\gO$-invariant, $\phi$ is homomorphism
from the $\gO$-module $\Sb^\ell$ into its contragredient module
$(\Sb^\ell)^*$. Both are simple $\gO$-modules and $\phi$ is non zero
since $C_\ell(\Sb^\ell, \Sb^\ell) = \Sb^0 = \KK$ by Theorem
\ref{3.3}(4), so $\phi$ is an isomorphism by Schur's Lemma and this 
proves that $\kappa$ is non degenerate.

Assume now that $\kappa'$ is an $\ad$-invariant bilinear form. Then
the map $F \mapsto \kappa'|_F$ where $\kappa'_F(G) = \kappa(F,G)$ is a
homomorphism of $\gO$-modules from $\Sb^\ell$ into $(\Sb^m)^*$. When
$\ell \neq m$, $\Sb^\ell$ and $(\Sb^m)^* \equiv \Sb^m$ are not
isomorphic, so $\kappa'(\Sb^\ell, \Sb^m) = \{ 0 \}$ by Schur's
Lemma. If $\ell = m$, $\Sb^\ell$ is a simple highest weight
$\gO$-module, so it is Schur irreducible, and it results that there is
on $\Sb^\ell$ only one invariant bilinear form, up to a scalar, and
consequently there exists $\alpha_\ell \in \KK$ such that
$\kappa'|_{\Sb^\ell \times \Sb^\ell} = \alpha_\ell \kappa|_{\Sb^\ell
  \times \Sb^\ell}$. Next we want to prove that $\alpha_\ell =
\alpha_{\ell+1}$ for all $\ell \geq 0$. First we prove that $\alpha_{
  2\ell -1} = \alpha_{2 \ell}$ for all $\ell \geq 1$. Note that
$\ad(p_1)(q_1^{2\ell}) = 2 \ell q_1^{2 \ell -1}$, then from
$\kappa'(\ad(p_1)(q_1^{2 \ell}), F) = - \kappa'( q_1^{2\ell} ,
\ad(p_1)(F))$, for all $F \in \Sb^{2\ell-1}$, we deduce that:
\[(\alpha_{ 2\ell-1} - \alpha_{2 \ell}) \kappa(\ad(p_1)(q_1^{2 \ell}),
F)= 0,\]

and since $\kappa$ is non degenerate on $\Sb^{2 \ell -1} \times \Sb^{2 \ell -1}$,
one can conclude that $\alpha_{ 2\ell -1} = \alpha_{2 \ell}$ for all
$\ell \geq 1$.

We then show that $\alpha_{ 2\ell} = \alpha_{2 \ell +1}$ for all $\ell
\geq 0$. One has $\kappa'(\ad(p_1^3)(q_1^{2 \ell}) , F) = - \kappa'(q_1^{2
  \ell} , \ad(p_1^3)(F))$, for all $F \in \Sb^{2\ell+1}$ and that
implies 
\[(\alpha_{ 2\ell+1} - \alpha_{2 \ell}) \kappa(\ad(p_1^3)(q_1^{2
  \ell}), F)= 0,\]

for all $F \in \Sb^{2\ell+1}$. But $\ad(p_1^3)(q_1^{2\ell})$ has
component $6 \ell p_1^2 q_1^{2 \ell -1}$ on $\Sb^{2 \ell +1}$ so
$\alpha_{ 2\ell} = \alpha_{2 \ell +1}$ as wanted.

Now, starting from $\ell =1$, we conclude that $\kappa' = \alpha_1 \kappa$
on $[\Wb, \Wb]$.
\end{proof}

\begin{cor}
$\Wb$ and $[\Wb, \Wb]$ are superquadratic Lie superalgebras.
\end{cor}

\begin{rem}
Consider the adjoint representation of $\g$ in $\Wb$. By Proposition
\ref{2.4}, it decomposes into isotypic components $\Wb = \bigoplus_{k
  \geq 0} A_k$, with $A_0 = \KK$, $A_k = \Sb^{2k-1} \oplus \Sb^{2k}$,
$k \geq 1$. Each $A_k$ has an explicit invariant supersymmetric
bilinear form, namely $\kappa|_{A_k \times A_k}$, so the
$\ad$-representation of $\g$ in $A_k$ is valued in
\[\osp (\Sb^{2k}, \Sb^{2k-1}) = \osp \left( \binom{2n+2k-1}{2k},
\binom{2n+2k-2}{2k-1} \right).\]

Also the $\ad$-representation of $\gO$ in $\Sb^k$ is orthogonal or
symplectic, according to the parity of $k$: $\Sb^{2k} $ is orthogonal,
$\Sb^{2k+1}$ is symplectic, and the corresponding bilinear form is
explicitly computable using the Moyal $\star$-product.
\end{rem}

What about the $\ad'$-representation of $\Wb$? In what follows, we
shall prove that it has also a non degenerate supersymmetric invariant
bilinear form. Actually, this bilinear form extends the one defined on
$\Sb^0 \oplus \Sb^1$ when embedding $\osp(1,2n)$ in $\Wb$ (see
Proposition \ref{2.3}), so it is directly related to orthosymplectic
supersymmetry.

\begin{thm} \label{4.5} 
Let $B(F,G) = (-1)^{fg+1} \kappa(F,G)$, for all $F \in \Sb^f$, $G
\in \Sb^g$. Then $B$ is a non degenerate supersymmetric bilinear
form on $\Wb$. Moreover, $B(\Sb^\ell, \Sb^m) = \{ 0 \}$, if $\ell
\neq m$, $B_{\Sb^\ell \times \Sb^\ell}$ is non degenerate and
$B$ is invariant under the $\ad'$-representation of $\Wb$. If
$B'$ is an $\ad'$-invariant bilinear form on $\Wb$, there exists
$\beta \in \KK$ such that $B' = \beta \ B$.
\end{thm}

\begin{proof}
it is easy to check that $B$ is supersymmetric. Let us prove that
$B$ is $\ad'$-invariant: consider $I = B(\ad'(A)(F), G) +
(-1)^{af} B(F, \ad'(A)(G))$ with $\deg_{\ZZ_2}(A) = a$,
$\deg_{\ZZ_2}(F) = f$ and $\deg_{\ZZ_2}(G) = g$. If $f + g + a =
\overline{1}$, then $I=0$. If $f + g + a = \overline{0}$, then
\begin{eqnarray*}
I &=& (-1)^{(f+a)g+1} \kappa(\ad'(A)(F), G) + (-1)^{af}
(-1)^{(f+a)g+1}\kappa(F, \ad'(A)(G)) \\ &=& (-1)^{f} \left( (-1)^{a+1}
\kappa(\ad'(A)(F), G) + (-1)^{af+1} \kappa(F, \ad'(A)(G)) \right)
\\ &=& (-1)^{f} \Str \left( (-1)^{a+1} (A\sta F - (-1)^{a(f+1)} F \sta
A) \sta G +
\right. \\ && \phantom{@@@@@@@@@@} \left. (-1)^{af+1} F \sta (A \sta G -
(-1)^{a(g+1)} G \sta A) \right) \\ &=& (-1)^{f+a+1} \Str \left( A \sta F
\sta G - (-1)^{a(f+g)} F \sta G \sta A \right) \\ &=& (-1)^{f+a+1} \Str (
     [A,F \sta G]) =0 
\end{eqnarray*} 

Since $\kappa$ is non degenerate, $B$ is non degenerate. Finally, in
order to prove the uniqueness of two $\ad'$-invariant bilinear forms
modulo $\KK$, one proceeds as in the proof of Theorem \ref{4.4}, so we
leave out the proof.
\end{proof}

\section{Renormalized supertrace and formal inverse Weyl
  transform} \label{Section05} 

In this Section, we assume that $\KK=\RR$ or $\CC$. Let $\Pc$ be the
algebra $\KK[x_1, \dots, x_n]$, $\Fc$ be the algebra $\KK[[X_1, \dots,
    X_n]]$, $V$ be the space $V = \spa\{X_1, \dots, X_n\}$ and $V^*$
be its dual, $V^* = \spa \{x_1, \dots, x_n\}$ with $\langle x_i, X_j
\rangle = \delta_{ij}$.

There is a one to one mapping (the Laplace transform) from $\Pc^*$
onto $\Fc$ defined by the duality $\langle x^I, X^J \rangle =
\delta_{IJ} I!$, where $x^I := x_1^{i_1} \dots x_n^{i_n}$, $X^J :=
X_1^{j_1} \dots X_n^{j_n}$, $I! := i_1! \dots i_n!$. So the spaces
$\Pc^*$ and $\Fc$ can be identified and as a consequence, the Dirac
distributions $\partial_v$, $v \in V$, $\partial_v(P) = P(v)$ become
formal power series $\exx^v$ so that:
\[ P(v) = \langle P \mid \exx^v \rangle, \ \forall \ P \in \Pc, v \in
V \ \text{ (Taylor's Formula). }\]

Also one has :
\[\left\langle \left. \dfrac{\partial^I P}{\partial x^I} \right| F
\right\rangle = \langle P \mid X^I F \rangle., \ \forall \ P \in \Pc,
F \in \Fc\]
 
and that means ${\phantom{\frac11}}^t \left(
\frac{\partial^I}{\partial x^I} \right)$ is the multiplication by
$X^I$ in $\Fc$.

The algebra $\Pc$ has a Hopf algebra structure with coproduct $\Delta(P)
:= P(x + x')$ if one identifies $\Pc \otimes \Pc = \KK [x_1, \dots,
  x_n, x_1', \dots, x_n' ]$, and antipode $\Sc(P)(x) =P(-x)$. 

Next we endow $\Pc$ with its natural topology, as defined in
\cite{BFGP} and $\Fc = \Pc^*$ with the strong dual topology which is
exactly the product topology $\Pi_{k \geq 0} \Fc^k$, where $\Fc^k$
denotes the set of homogeneous polynomials of degree $k$. Then the
transposition map induces on $\Fc$ a topological Hopf algebra
structure, which is exactly the usual structure (see \cite{BFGP}) with
the identification $\Fc \widehat{\otimes} \Fc = \KK[[X_1, \dots,
    X_n,X_1', \dots, X_n']]$ (where $\widehat{\otimes}$ is the
projective tensor product, see \cite{Treves}).

Any linear operator $T \colon \Pc \to \Pc$ is continuous for the
natural topology (\cite{BFGP}) and ${}^t T \colon \Fc \to \Fc$ is
continuous. Denote by $\Lc(\Pc)$ the space of linear operators of
$\Pc$, and by $\Lc(\Fc)$ the space of continuous linear operators of
$\Fc$. Then one has:
\[ \Lc(\Pc) = \Pc^* \widehat{\otimes} \Pc = \Pc^*
\widehat{\otimes} \Pc^{**} = \Lc(\Fc) \quad \text{ (see
  \cite{Treves}).}\]

Let us quickly explain how it works: given $T \in \Lc(\Pc)$, then $T =
\sum_K \frac{1}{K!} P_K \otimes X^K$ with $P_K = T(x^K)$ and one has
$T(P) = \sum_K \frac{1}{K!} P_K \langle X^K \mid P \rangle$ for all $P
\in \Pc$. Now one has ${}^tT = \sum_K \frac{1}{K!} P_K \otimes X^K$
and 
\begin{equation}\label{5.1}
{}^tT(F) = \sum_K \frac{1}{K!} \langle P_k \mid F \rangle X^K \text{
  for all } F \in \Fc.
\end{equation} 

Given a polynomial function $d \colon V \to \Fc$
defined as $d (v) = \sum_K D_K(v) X^K$, $D_K \in \Pc$, there is an
associated operator $D \colon \Pc \to \Pc$ defined by 
\[D = \sum_K D_K(v) \frac{\partial^K}{\partial x^K} \]

and one has:
\[ D(P)(v) = \langle P \mid d (v) \exx^v \rangle, \ \forall  \ P \in
\Pc, v \in V.\]

Since $\spa \{ \exx^v \mid v \in V \}$ is dense in $\Fc = \Pc^*$ by
the Hahn-Banach Theorem, we deduce:
\begin{equation}\label{5.4} 
\langle D(P) \mid F \rangle = \langle P \mid m \circ (d \otimes \Id)
\circ \Delta (F) \rangle,
\end{equation}

noticing that $\Delta(\exx^v) = \exx^v \otimes \exx^v$ and extending
$d$ to $\Fc$ by $d(F) = \sum_K \langle D_k \mid F \rangle \ X^k$, $F
\in \Fc$ (so $d(\exx^v) = d(v)$).

The operator $D$ is what we shall call a {\em differential operator }
of $\Pc$. We shall say that $D$ is {\em differential operator for
  finite order} if $d(V) \subset \Fc_k = \{ F \in \Fc \mid \deg(F)
\leq k \}$.

A fundamental property of $\Pc$ is established by:

\begin{lem} \label{5.6} (\cite{Pinczon1})
Any linear operator of $\Pc$ is a differential operator.
\end{lem}

\begin{proof}
Given $T \in \Lc(\Pc)$, then ${}^t T \in \Lc(\Fc)$, so we have to find
a polynomial map from $V$ to $\Fc$ satisfying ${}^tT = m \circ (d
\otimes \Id) \circ \Delta$ (due to (\ref{5.4})) From the density of
$\spa \{ \exx^v \mid v \in V \}$, it is enough to prove the last
identity on this set. Let $d = \sum_k D_K X^K$, $D_K \in \Pc$, then
${}^t T (\exx^v) = (m \circ (d \otimes \Id) \circ \Delta) (\exx^v) =
\sum_K \langle D_k \mid \exx^v \rangle X^K \exx^v$ gives $\sum_K
\langle D_K \mid \exx^v \rangle X^K = {}^t T ( \exx^v) \exx^{-v}$,
hence $d(v) = {}^t T(\exx^v) \exx^{-v}$.
\end{proof}

Let us give an explicit formula: starting from $T \in \Lc(\Pc)$, $T =
\sum_I \frac{1}{I!} P_I \otimes X^I$ with $P_I \in \Pc$. Then ${}^t
T(\exx^v) = \sum_I \frac{1}{I!} P_I(v) X^I$ by (\ref{5.1}) so if $v =
x_1 X_1 + \dots x_n X_n$ and $|I| = i_1 + \dots + i_n$,
\begin{eqnarray*}
&&{}^t T(\exx^v) \exx^{-v} = \\ &=&\left( \sum_I \frac{1}{I!} P_I(v)
  X^I \right) \left( \sum_j (-1)^j \dfrac{v^j}{j!} \right) \\ &=&
  \left( \sum_I \frac{1}{I!} P_I(v) X^I \right) \left( \sum_j
  \dfrac{(-1)^j}{j!}  \sum_{i_1 + \dots + i_n=j} \dfrac{(i_1 + \dots +
    i_n)!}{i_1! \dots i_n!} x_1^{i_1} \dots x_n^{i_n} X_1^{i_1} \dots
  X_n^{i_n} \right)\\ &=& \left( \sum_I \frac{1}{I!} P_I(v) X^I
  \right) \left( \sum_I \frac{(-1)^{|I|}}{I!} x^I X^I \right) \\ &=&
  \sum_N \frac{1}{N!}  \left( \sum_{R+S=N} (-1)^{|S|} \dfrac{N!}{R!S!}
  P_R(v) x^S \right) X^N.
\end{eqnarray*}

Finally 
\[T =\sum_N \frac{1}{N!}  \left( \sum_{R+S=N} (-1)^{|S|}
\dfrac{N!}{R!S!} P_R x^S \right) \dfrac{\partial^N}{\partial x^N}.\]

Recall that $P_R = T(x^R)$, therefore:
\[T = \sum_N \frac{1}{N!} \left( \sum_{R+S=N} (-1)^{|S|}
\dfrac{N!}{R!S!} T(x^R) x^S \right) \dfrac{\partial^N}{\partial
  x^N}.\] 

But $\Delta(x^N) = \sum_{R+S=N} \dfrac{N!}{R!S!} x^R \otimes x^S$, so
this last formula can be written:

\begin{thm} \label{5.7.2}
\[ T =  \sum_N \frac{1}{N!} \left( m \circ (T \otimes \Sc) \circ \Delta(x^N)
\right) \ \dfrac{\partial^N}{\partial x^N}. \]
\end{thm}

where $\Sc$ is the antipode of $\Pc$. 

For instance, when $n =1$, consider the operator $T$ defined as
$T(x^i) = x^j$ for fixed $i$ and $j$, and $0$ otherwise. Then;
\[T =  \frac{x^j}{i!} \sum_{\ell \geq 0} (-1)^{\ell}
\dfrac{x^\ell}{\ell!} \dfrac{d^{i+\ell}}{d x^{i+\ell}}.\]

Using the notation $P^I = p_1^{i_1}\dots p_n^{i_n}$ and $Q^J =
q_1^{j_1}\dots q_n^{j_n}$, consider the (formal) completion
$\overline{\Wb} = \KK [Q][[P]]$ of the $\star$-algebra $\Wb$
($\overline{\Wb}$ is an algebra with the Moyal
$\star$-product). Elements of $\overline{\Wb}$ are formal power series
$\widetilde{F} = \sum_I \alpha_I(Q) \sta P^I$ with $\alpha_i \in
\KK[Q]$. Define a map $\Wc \colon \overline{\Wb} \to \Lc(\Pc)$ as:
\[\Wc (\widetilde{F}) = \sum_I \alpha_i(x)
\ \dfrac{\partial^I}{\partial x^I}.\] 

Remark that $\Wc$ is simply the extension to $\overline{\Wb}$ of the
natural $\Wb$-module structure of $\Pc$ defined by $p_i \cdot P =
\dfrac{\partial P}{\partial x_i}$, $q_i \cdot P = x_i P$, $i = 1,
\dots, n$, hence $\Wc$ is an algebra homomorphism.

Lemma \ref{5.6} can now be reinterpreted as:

\begin{prop} \label{5.7}
The map $\Wc$ is an isomorphism of algebras.
\end{prop}

One should be careful that the domain of $\Wc$, i.e. $\overline{\Wb}$,
is not at all identical to the formal completion $\overline{S} =
\KK[Q][[P]]$ endowed with an Abelian product and with elements
$\widehat{F} = \sum_I \alpha_i(Q) \ P^I$.

Given $\widetilde{F} = \sum_I \alpha_i(Q) \sta P^I \in
\overline{\Wb}$, we can try to define its supertrace using the
supertrace of $\Wb$ that we will denote by $\Str_\Wb$. A natural
candidate would be $\Strwb(\widetilde{F}) := \sum_I \Strw(\alpha_i(Q)
\sta P^I)$ but it is clear that this series happens to diverge. So
$\Strwb$ has a domain denoted by $\Dom(\Strwb)$. Evidently, $\Wb
\subset \Dom(\Strwb)$. On the other hand, $\Lc(\Pc)$ is a
$\ZZ_2$-graded algebra since $\Pc$ is $\ZZ_2$-graded, so we can define
a supertrace at least on the ideal $\Lc_f(\Pc)$ of finite rank
operators. Note that $\Lc_f(\Pc) \cap \Wb = \{ 0 \}$ since $\Wb$ is a
simple algebra. So we now have two supertraces, living apparently on
different domains, and we wish to compare these supertraces. This is
done by the following Theorem:

\begin{thm}\label{5.9} \hfill

One has $\Lc_f(\Pc) \subset \Dom(\Strwb)$ and if $\widetilde{F} \in
  \Lc_f(\Pc)$,
\[ \Str(\widetilde{F}) = \frac{1}{2^n} \ \Strwb(\widetilde{F}).\]
\end{thm}

To prove Theorem \ref{5.9}, we need precise formulas for $\Strw$:

\begin{prop}  \label{5.9.1} \hfill

\begin{enumerate}
\item Using the natural isomorphism $\Wb(n) = \Wb(1) \otimes \dots
  \otimes \Wb(1)$, one has
\[ \Strw(F_1 \otimes \dots \otimes F_n) = \Strw(F_1) \dots \Strw(F_n).\]

\item One has $\Strw(P^I \sta Q^J) = \delta_{IJ} \dfrac{I!}{2^{|I|}}$.

\end{enumerate}
\end{prop}

\begin{proof}\hfill

\begin{enumerate}
\item Recall that the isomorphism $\Wb(n) = \Wb(1) \otimes \dots
  \otimes \Wb(1)$ is defined by $F_1 \otimes \dots \otimes F_n = F_1
  \dots F_n$. Then 
\begin{eqnarray*}
\Strw(F_1 \otimes \dots \otimes F_n) &=& \Strw(F_1 \dots F_n) = (F_1
\dots F_n) (0) = F_1(0) \dots F_n(0) \\ &=& \Strw(F_1) \dots
\Strw(F_n).
\end{eqnarray*}

\item We start with the case $n =1$. Then $\Strw(p_1^{i_1} \sta
  q_1^{j_1}) = \kappa(p_1^{i_1} , q_1^{j_1}) = 0$ if $i_1 \neq j_1$ by
  Proposition (\ref{4.4}). Furthermore, $\Strw(p_1^{i_1} \sta
  q_1^{i_1}) = C_{i_1}(p_1^{i_1} , q_1^{i_1}) =
  \dfrac{i_1!}{2^{i_1}}$. Compute now:
\begin{eqnarray*}
\Strw(p_1^{i_1} \dots p_n^{i_n} \sta q_1^{j_1} \dots q_n^{j_n} ) &=&
\Strw(p_1^{i_1}\sta q_1^{j_1} \otimes p_2^{i_2}\sta q_2^{j_2} \otimes
\dots \otimes p_n^{i_n} \sta q_n^{j_n} ) \\ &=& \Strw(p_1^{i_1}\sta
q_1^{j_1}) \dots \Strw(p_n^{i_n} \sta q_n^{j_n} ) =0,
\end{eqnarray*}

if $I \neq J$, and $\dfrac{I!}{2^{|I|}}$ if $I = J$.
\end{enumerate}
\end{proof}

As a practical case of Theorem \ref{5.9}, we prove:

\begin{prop} \label{5.9.2}
Let $T$ be the (finite rank) operator defined as $T(x^I) = x^J$ for
fixed $I$, $J$ and $0$ otherwise. Then
\[T = \frac{x^J}{I!} \sum_S (-1)^S \frac{x^S}{S!}
\frac{\partial^{I+S}}{\partial x^{I+S}}, \]

and 
\[ \Strwb(T) = \begin{cases} 2^n \, \Str(T) = 2^n (-1)^{|I|}, \text{ if
  } I =J\\0 \text{ otherwise } \end{cases}.\]

Thus $T \in \Dom(\Strwb)$ and its usual supertrace is, up to a factor
$\frac {1}{2^n}$, its $\overline{\Wb}$-supertrace.
\end{prop}

\begin{proof}
The formula for $T$ is obtained by applying Theorem \ref{5.7.2}:
\[T = \sum_S (-1)^{|S|} \dfrac{Q^{J+S}}{I!S!} \sta P^{I+S}.\]

Then $\Strwb(T) = \sum_S \dfrac{(-1)^{|S|}}{I!S!} \Strw( Q^{J+S} \sta
P^{I+S}) =0.$, if $J \neq I$, by Proposition \ref{5.9.1}. When $I = J$,
$\Strwb( Q^{I+S} \sta P^{I+S}) = (-1)^{|I+S|}
\dfrac{(I+S)!}{2^{|I+S|}}$, so
\[\Strwb(T) = \dfrac{(-1)^{|I|}}{2^{|I|}} \sum_S  \dfrac{(I+S)!}{I!S!}
\dfrac{1}{2^{|S|}}.\]

If $|\tau| < 1$, one has $\dfrac{1}{(1-
    \tau)^{i+1}} = \sum_s \binom{i+s}{s} \tau^s$, therefore: 
\[\dfrac{1}{(1-\tau_1)^{i_1+1}\dots (1-\tau_n)^{i_n+1}} = \sum_{s_1,
  \dots,s_n} \binom{i_1+s_1}{s_1} \dots \binom{i_n+s_n}{s_n}
\tau_1^{s_1} \dots \tau_n^{s_n},\]

if $|\tau_i| <1$, for all $i$. It results that:
\[\sum_S  \dfrac{(I+S)!}{I!S!} \dfrac{1}{2^{|S|}} =
\dfrac{1}{\left( 1-\frac12 \right)^{i_1+1}\dots \left( 1-\frac12
  \right)^{i_n+1}} = 2^{|I|+n}.\]

At last, we obtain
\[\Strwb(T) = (-1)^{|I|} 2^n  =  \Str(T) 2^n.\]
\end{proof}

\begin{proof} (Proof of Theorem \ref{5.9}) \hfill

This proof is completely similar to the previous one and for this
reason, we omit it (note that one can restrict to $T = \varphi \otimes
x^k$, $\varphi \in \Pc^*$ since $\Lc_f(\Pc)= \Pc^* \otimes \Pc$). 
\end{proof}

\begin{rem}
Let us give some more interpretation about supertraces. We want to
show how one can define a renormalized supertrace from Theorem
\ref{5.9}. We begin with $\Lc_f(\Pc)$ and its natural supertrace
$\Str$. Notice that this supertrace is of intrinsic nature, since
defined by $\Str(\varphi \otimes P) = (-1)^{\deg_{\ZZ_2}(\varphi)
  \deg_{\ZZ_2}(P)} \langle \varphi \mid P \rangle$, $\varphi \in
\Pc^*$, $P \in \Pc$. Secondly, we have the supertrace $\Strwb$ defined
on its domain $\Dom(\Strwb)$ which contains $\Wb$ and $\Lc_f(\Pc)$. On
$\Lc_f(\Pc)$, one has 
\[ \Str(T) =\frac{1}{2^n} \Strwb(T)\]

by Theorem \ref{5.9}. So we can extend $\Str$ to $\Dom(\Strwb)$ and
define a {\em renormalized supertrace}, denoted by $\RStr$, as:
\[ \RStr(T) = \frac{1}{2^n} \Strwb (T).\]

This extension is indeed a renormalized extension of $\Str$: for
instance, with the usual definition of the supertrace: $\Str(\Id) =
\infty - \infty$, a rather bad result, but with the renormalization:
\[ \RStr(\Id) = \frac{1}{2^n}.\]

Notice that $n = \dim(V)$ is the dimension of the underlying variety.
\end{rem}

Next, we will clarify what we mean by a {\em formal inverse Weyl
  transform}. Recall that we can identify $\overline{\Wb} = \Lc(\Pc)$
thanks to Proposition \ref{5.7}. So, given $T \in \Lc(\Pc)$, one can
write:
\[T=\sum_I \alpha_I(Q) \sta P^I \in \overline{\Wb}.\]

Consider the Abelian algebra $\KK[[P,Q]]$, and denote by
$\KK[[P,Q]]^k$ the space of homogeneous polynomials of degree $k$. We
endow $\KK[[P,Q]]$ with the product topology $\KK[[P,Q]] =
\Pi_{k\geq0} \KK[[P,Q]]^k$ which is Fréchet. Now $\sum_I \alpha_I(Q)
\sta P^I \in \overline{\Wb}$ is a series (of polynomials) in
$\KK[[P,Q]]$, it converges in $\KK[[P,Q]]$ if and only if for any $k$,
the series obtained by taking the $k$th components converge in
$\KK[[P,Q]]^k$.

\begin{defn} \label{5.8}
When $\sum_I \alpha_I(Q) \sta P^I$ converges in $\KK[[P,Q]]$, we say
that $T$ has a {\em formal inverse Weyl transform} denoted by $\IW(T)$
and defined as the sum of the series.
\end{defn}

Remark that taking $k=0$, we obtain that if $T$ has a formal inverse
Weyl transform, the series 
\[\sum_I \left( \alpha_I(Q) \sta P^I \right) (0)
= \sum_I \Strw \left( \alpha_I(Q) \sta P^I \right) = \Strwb(T)\]

must converge, so:

\begin{prop}
The existence of $\RStr(T)$ is a necessary condition for the existence
of the formal inverse Weyl transform of $T$. When $\IW(T)$ exists, one
has:
\[ \RStr(T) = \dfrac{1}{2^n} \IW(T)(0).\]
\end{prop}

As an operator from $\Lc(\Pc)$ into $\KK[[P,Q]]$, the formal inverse
Weyl transform has a domain $\Dc \subsetneq \Lc(\Pc)$ containing $\Wb$
and naturally, it would be nice to have a characterization of
$\Dc$. This will be done elsewhere, rather let us develop here some
examples. We take $n =1$ and consider the elementary operators
$E_{ij}$ of $\Pc$ defined by $E_{ij} (x^k) = \delta_{jk} x^i$. One has
$E_{ij} = \frac{1}{j!} \sum_{\ell \geq 0} \frac{(-1)^\ell}{\ell!}
q^{\ell +i} \sta p^{\ell+j}$ and some computation shows that

\begin{lem}
The formal inverse Weyl transform of $E_{ij}$ is
\[ \IW(E_{ij}) = \begin{cases} (-1)^j 2^{i-j+1} L_j^{(i-j)} (4pq)
  \exx^{-2pq} q^{i-j}, \text{ if } j \leq i \\ (-1)^i 2^{j-i+1}
  \frac{i!}{j!} L_i^{(j-i)} (4pq) \exx^{-2pq} p^{j-i}, \text{ if } j
  \geq i \end{cases}.\]
\end{lem}

where $L_\beta^{(\alpha)}$ is the Laguerre polynomial (see
\cite{Szego}).

Consider now the operators $S_\lambda$ defined as $S_\lambda (x^k) =
\lambda^k x^k$, $\lambda \in \KK$. Then by Theorem \ref{5.7.2}, one has:
\[S_\lambda = \sum_{\ell \geq 0} \frac{(-1)^\ell}{\ell !} (1 -
\lambda)^\ell q^\ell \sta p^\ell.\]

By (\ref{1.9.2}) in Section \ref{Section01}, $q^\ell \sta p^\ell =
(-1)^\ell \frac{\ell !}{2^\ell} L_\ell (2pq)$, so using the generating
function of Laguerre polynomials $L_\beta^{(0)} := L_\beta$, one
finds:

\begin{lem} \label{5.11.3}
The formal inverse Weyl transform of $S_\lambda$ is 
\[ \IW(S_\lambda)=\dfrac{2}{1+\lambda} \exp \left( 2 \dfrac{\lambda - 1}{\lambda+1}
pq \right), \text{ if } |1-\lambda|< 2.\]
\end{lem}

Notice that $S_\lambda$ has a renormalized supertrace given by:
\[ \RStr(S_\lambda) = \dfrac{1}{\lambda+1}.\]

When $|\lambda| < 1$, this is exactly $\Str(S_\lambda)$. 

Here are some cases of interest:

First, take $\lambda = \exx^\tau$, then $S_\lambda = \exp \left( \tau
x \dfrac{d}{dx} \right) = \exp (\tau q \sta p)$ and one finds
that its formal inverse Weyl transform is:
\[ \IW \left( \exp \left( \tau
x \dfrac{d}{dx} \right) \right) = \frac{\exx^{-\frac{\tau}{2}}}{\cosh
  \left(\frac{\tau}{2} \right)} \exp \left( 2pq \tanh \left(
\frac{\tau}{2} \right) \right) , \text{ if } |1-\exx^\tau|< 2.\]

Take $\tau = i \theta$, $\theta \in \RR$ with $\theta \neq (2s+1)
\pi$, $s \in \ZZ$. Then $S_{\exx^{i \theta}}= \exp ( i \theta q \sta
p)$ has formal inverse Weyl transform:
\[ \IW \left( \exp \left( i \theta
x \dfrac{d}{dx} \right) \right) = \dfrac{\exx^{-i
    \frac{\theta}{2}}}{\cos \left(\frac{\theta}{2} \right)} \exp
\left( 2 i pq \tan \left( \dfrac{\theta}{2} \right) \right) ,\]

and an interesting case is $S_i$ with formal inverse Weyl transform
$(1-i) \exp(2i pq)$ (see \cite{BFFLS}, where similar formulas are
obtained, in the context of quantization of the harmonic
oscillator). Note that the value $\lambda = -1$ is critical: indeed
$S_{-1} =P$ is the parity operator, which has divergent renormalized
supertrace and cannot have a formal inverse Weyl transform. So the
estimation in (\ref{5.11.3}) is the best one.

Let us now justify our construction somewhat. There is a natural
question : quantum mechanics is built from operators in a Hilbert
space, so where is the Hilbert space in that picture? Here is the
answer: let $\Hc$ be the Hilbert space of entire functions $f(x)$, $x
= ai+b$ such that $\int \exx^{-|x|^2} |f(x)|^2 \ da \ db < \infty$
with scalar product $\langle f \mid g \rangle = \frac{1}{\pi} \int
\exx^{-|x|^2} f(x) \overline{g(x)} \ da \ db$. Then $\Pc \subset \Hc$,
as a dense subspace and we can apply our algebraic formalism as
follows: given $T$ an operator of $\Hc$, whose domain contains $\Pc$,
assume that $T(\Pc) \subset \Pc$ and denote by $T_\Pc$ the restriction
to $\Pc$; to have one to one $T \mapsto T_\Pc$, we have to assume
more, for instance either $T$ bounded, or $T_\Pc$ essentially
self-adjoint, this is assumed in the foregoing. We can now define the
formal inverse Weyl transform of $T$ to be the formal inverse Weyl
transform of $T_\Pc$ and the renormalized supertrace of $T$ to be the
renormalized supertrace of $T_\Pc$. For instance, consider $S_i$ which
is a unitary operator of $\Hc$, it has renormalized supertrace
$\RStr(S_i) = \frac12 (1-i)$ and formal inverse Weyl transform $(1-i)
\exp(2i pq)$. We are actually working in the so called ``coherent
states formalism'', and we can easily translate in terms of the usual
``metaplectic formalism''. Define an operator $H \colon \Hc \to
\Lc^2(\RR)$ by
\[H(f)_t = \dfrac{\exx^{-\frac{t^2}{2}}}{\pi^{\frac34}} \int f(\sqrt{2}
(t + i \xi)) \exx^{-\xi^2} d \xi,\]

let $\phi_n$ be the orthonormal basis of Hermite functions in
$\Lc^2(\RR)$, and $Z_n = \frac{x^n}{\sqrt{n!}}$ the orthonormal basis
of $\Hc$, then one has $H(Z_n) = \phi_n$, $\forall n$, so $H$ is a
unitary isomorphism from $\Hc$ to $\Lc^2(\RR)$. The operator $H$ maps
the operator $x$ (resp. $\frac{d}{dx}$) or $\Hc$ onto the operator
$\frac{1}{\sqrt{2}} \left( t - \frac{d}{dt} \right)$
(resp. $\frac{1}{\sqrt{2}} \left( t + \frac{d}{dt} \right)$) and we
recover the usual metaplectic formalism (see \cite{Pi-Le} for
details). now we can explain our interest for $S_i$: indeed $H$ maps
$S_i$ onto the Fourier transform $\mathfrak{F}$ of $\Lc^2(\RR)$ and
since ``coherent states formalism'' and ``metaplectic formalism'' are
completely equivalent by $H$, we have a formal inverse Weyl transform
of $\mathfrak{F}$, namely the formal inverse Weyl transform of $S_i$,
i.e.
\[ \IW(\mathfrak{F}) = (1-i) \exp(2ipq)\]

(attention to the interpretation of $p$ and $q$ in terms of $t$), and
also a renormalized supertrace 
\[\RStr(\mathfrak{F}) = \frac12 (1-i).\]

So we think that our formalism might be of interest. Next step will be
a Wigner's type formula: this is already done and will be explained in
a a subsequent paper.

\end{document}